\input amstex
\input epsf
\documentstyle{amsppt}
\magnification=1200
%\hcorrection{0.5in}
%\vcorrection{-0.6in}
%\hcorrection{0.7in}
%\vcorrection{0.7in

%%\predefine\Vol{\vol}
%\define\volume{\operatorname{vol}}
%\define\svol{\underline{\operatorname{vol}}}
%\define\svolball#1#2{{\svol(B_{#2}^{#1})}}
%\define\svolann#1#2{{\svol(A_{#2}^{#1})}}
%\define\sball#1#2{{B_{#2}^{#1}}}
%\define\svolsp#1#2{{\svol(\partial B_{#2}^{#1})}}
\define\volume{\operatorname{vol}}
\define\op#1{\operatorname{#1}}
%the following is the standard space volume of ball, annulus, sphere, please use them instead of typing concrete symbols.
\define\svolball#1#2{{\volume(\underline B_{#2}^{#1})}}
%standard volume of ball of radius #2 in constant curvature #1
\define\svolann#1#2{{\volume(\underline A_{#2}^{#1})}}
%standard volume of annulus of radius #2 in constant curvature #1
\define\sball#1#2{{\underline B_{#2}^{#1}}}
%standard ball of radius #2 in constant curvature #1
\define\svolsp#1#2{{\volume(\partial \underline B_{#2}^{#1})}}
%standard volume of sphere of radius #2 in constant curvature #1

\documentstyle{amsppt}
\magnification=1200 \hcorrection{0.2in} \vcorrection{-0.1in}
\NoRunningHeads \NoBlackBoxes \topmatter
\title The Soul Conjecture in Alexandrov Geometry in dimension 4
\endtitle
\author Xiaochun Rong \footnote{Supported partially by
NSF Grant DMS 1106517 and a research found from Capital normal university.
\hfill{$\,$}} \& Yusheng Wang \footnote{Supported partially by
NFSC 11471039.\hfill{$\,$}}
\endauthor
\abstract   In this paper, we prove the Soul Conjecture in
Alexandrov geometry in dimension $4$, i.e. if $X$ is a complete non-compact $4$-dimensional
Alexandrov space of non-negative curvature and positive curvature around one point, then a soul of $X$ is a point.
\endabstract

\address Mathematics Department, Capital Normal University,
Beijing, P.R.C. \newline.\hskip3mm
Mathematics Department, Rutgers University, New Brunswick,
NJ 08903, U.S.A
\endaddress

\address School of Mathematical Sciences (and Lab. math. Com.
Sys.), Beijing Normal University, Beijing, 100875 P.R.C. ({\rm E-mail: wyusheng\@bnu.edu.cn})
\endaddress
\endtopmatter

\document

\vskip2mm

%%%%%%%%%%%%%%%%%%%%%%%%%%%%%%%%%%%%%% section 0 Introduction %%%%%%%%%%%%%%%%%%%%%%%%%%

\head 0. Introduction
\endhead

\vskip4mm

The purpose of this paper is to prove the Soul Conjecture in
Alexandrov geometry ([Per1]) in dimension $4$.

Let's start with some background. In Riemannian geometry,
the classical Soul Theorem of Cheeger-Gromoll is ([CG], cf. [GM1]):

\proclaim{Theorem 0.1} Let $M$ be a complete non-compact Riemannian
manifold of non-negative sectional curvature. Then there is a compact totally convex
submanifold, $S$ (called a soul of $M$), such that $M$ is diffeomorphic
to the normal bundle of $S$.
\endproclaim

The following is the Soul Conjecture of Cheeger-Gromoll ([CG]), proved more than 20 years later by Perel'man ([Per1]):

\proclaim{Theorem 0.2} Let the assumptions be as in Theorem 0.1.
If there is open set $U\subset M$ on which sectional curvature $\op{sec}_U>0$, then $S$ is a point.
\endproclaim

Perel'man showed that if $S$ is not a point, then at every $x\in M$, there is a plane of
zero curvature. His approach relies on the existence of a distance non-increasing map,
$\pi: M\to S$, which can be taken as the Sharafutdinov retraction ([Sh]),
$\pi(x)=\lim_{t\to \infty}\phi_t(x)$, where $\phi_t$ is the Sharafutdinov flow determined by
a Busemann function $f$ on $M$ (see Section 1 for definitions). Precisely, Perel'man established
the following Flat
Strip property: for $x\in M\setminus S$, the distance $|x\pi(x)|=|xS|$, and a minimal geodesic from $x$ to $\pi(x)$ and any geodesic in $S$ at $\pi(x)$ bound
an isometrically embedded flat strip. At the same time, Perel'man proved that $\pi$ is
a $C^1$-Riemannian submersion (cf. [CS], [Wi]).

An Alexandrov space $X$ with curvature bounded below by $\kappa$ is a length metric space on which the Toponogov comparison holds with respect to a complete surface of constant curvature $\kappa$.
Since the seminal paper [BGP], there have been many works extending basic results in Riemannian geometry whose proofs rely on the Toponogov comparison to Alexandrov geometry; including Theorem 0.1
that any complete non-compact Alexandrov space of non-negative curvature homotopically retracts to
a soul $S$, a compact convex subset without boundary, via the Sharafutdinov retraction.

In [Per1], Perel'man asked whether Theorem 0.2 holds in Alexandrov geometry, which
we formulate as

\example{Conjecture 0.3} Let $X$ be a complete non-compact Alexandrov space of non-negative curvature
with a soul $S$. If there is an open subset $U$ such that curvature $\op{cur}|_U>0$, then $S$ is a point.
\endexample

An obstacle is that the Flat Strip property does not hold at every point (cf. [Li]). However,
Conjecture 0.3  will follow, if one can show that the Sharafutdinov retraction, $\pi: X\to S$, is a submetry,
i.e., for any $x\in X$ and any small $r>0$, $\pi(B_r(x))=B_r(\pi(x))$. This is  because $\op{cur}_U>0$
forces $\phi_t$ to strictly shrink the size of any metric ball in $U$.

The following cases have been known that the Sharafutdinov retraction, $\pi: X\to S$, is a submetry:

\noindent (0.4.1)  ([SY]) $\dim(S)=1$, or $\dim(S)=\dim(X)-1$.

\noindent (0.4.2) ([Li]) $\dim(S)=\dim(X)-2$ and $X$ is topologically nice, i.e. all iterated spaces
of directions at each point are homeomorphic to spheres ([Ka]).

Note that (0.4.1) implies Conjecture 0.3 in dimension $3$ ([SY]), and (0.4.2) implies a case of Conjecture 0.3 in dimension $4$ that $X$ is a topological manifold, because in dimension $4$,
$X$ is topological manifold iff $X$ is a topologically nice ([Li]); which does not hold in dimension $\ge 5$ ([Ka]).

The main results in this paper imply Conjecture 0.3 in dimension $4$.

\proclaim{Theorem A} Let $X$ be a complete non-compact Alexandrov $4$-space with curvature $\op{cur}\geq0$ and an empty boundary, and let $S$ be a soul of $X$. Then

\noindent {\rm (A1)} The Sharafutdinov retraction, $\pi: X\to S$, is a submetry.

\noindent {\rm (A2)} If there is an open set $U$ of $X$ such that $\op{cur}|_U>0$,
then $S$ is a point.
\endproclaim

Note that without loss of generality, we can assume that $X$ is simply connected
(cf. [Li]). In view of the above discussion, (A1) implies (A2),
and by (0.4.1) it remains to check the case of (A1): $S$ is homeomorphic to a sphere.

A common starting point, in verifying that $\pi: X\to S$ is a submetry, is the partial
Flat Strip property by Shioya-Yamaguchi
[SY]: for $p\in S$ and $\Omega_c\triangleq f^{-1}([c,c_0])$ with $c_0=\max f$, a minimal geodesic from $p$ to $\partial \Omega_c$
and a non-trivial minimal geodesic in $S$ at $p$ bound an isometrically embedded flat strip.
Let $\Uparrow_p^{\partial \Omega_c}$ denote the set of the directions of minimal geodesics from $p$ to $\partial \Omega_c$.
Note that $v\in \,\Uparrow_p^{\partial \Omega_c}$ actually points to
a ray at $p$, $\gamma_{v,p}$, that realizes the
distance to $\partial \Omega_c$ for all $c$ ([Li], see Section 1 below). Let
$$\Cal F=\bigcup_{p\in S, v\in \, \Uparrow_p^{\partial \Omega_c}}\gamma_{v,p}=\bigcup_{v\in \, \Uparrow_{p_0}^{\partial \Omega_c}}F_v,$$
where $F_v$ is the union of the  rays which are `parallel' to $\gamma_{v,p_0}$ through a piecewise flat strips over $S$ (see Section 3),
and $p_0\in S$ is a Riemannian point, i.e. the space of directions at $p_0$, $\Sigma_{p_0}S$, is isometric to the unit sphere.

By definition, points in $\Cal F$ satisfy the Flat Strip property,
and thus $\pi|_{\Cal F}: \Cal F\to S$ is a submetry. Hence throughout the rest of the paper, we will always assume that $\Cal F\subsetneq X$, and thus $\dim(S)<n-1$ (see (0.4.1)).

For $\dim(S)=n-2$, the main discovery in [Li] is that if $X$ is topologically nice, then $F_v$ is convex and splits, $F_v\cong S\times \Bbb R_+$ (here $S$ being simply connected is required).
Let $y\in \Cal F$ such that $|xy|=|x\Cal F|$. Because $\Cal F$ or $\partial \Cal F$ is a union of at most countably many $F_v$'s, $y\in F_v$ for some $v$. Observe that the distance function,
$d(\cdot, \Cal F): X\to \Bbb R_+$, is concave at $x$ if $y\in F_v\setminus S$. If $y\in S$, then the concavity for $d(\cdot, \Cal F)$ at $x$ can be easily checked \footnote{In [Li], Li actually
uses $d(\cdot, \bigcup  F_{v_i})$, where at most three $F_{v_i}$'s are suitably selected.}. When
$y\in F_v\setminus S$, because $F_v\setminus S$ locally
divides $X$ into two components, the gradient flow of $d(\cdot, \Cal F)$
extends to a distance non-increasing flow, $\Psi_t:F_v\setminus S\to X$, such that $\Psi_{|x\Cal F|}(y)=x$.

The significance of $F_v\cong S\times \Bbb R_+$ is a natural isometry: $\varphi_c: S\to \tilde S_c=F_v\cap \partial \Omega_c\cong S\times \{c\}$, $c=f(y)$.
Because the one-parameter family of composition maps, $\varphi_c\circ
\pi\circ \Psi_t: \tilde S_c\to \tilde S_c$, is onto and distance non-increasing,
$\varphi_c\circ\pi\circ \Psi_t$ is an isometry. Consequently,
$\Psi_t|_{\tilde S_c}: \tilde S_c\to \Psi_t(\tilde S_c)$ is an isometry,
and $\pi|_{\Psi_t(\tilde S_c)}: \Psi_t(\tilde S_c)\to S$ is a submetry. Hence, $\pi$ is a submetry at $x$ with $y\not\in S$, and $\pi$ is a submetry because the set of such $x$ is dense.

We point it out that the condition, `topologically nice',  is essential for [Li]; in general $F_v$ does not split, nor is it possible to have a distance non-increasing map, $S\to \tilde S_c$. Here is a simple example.

\example{Example 0.5} Let $M=T\Bbb S^2$ denote the tangent bundle of the unit sphere. Then the canonical metric on $M$ has non-negative sectional curvature. Let $D_{2m}\subset
SO(3)$ be the Dihedral group. Then the differential of
the $D_{2m}$-action defines an isometric $D_{2m}$-action on $M$ which preserves
$\Bbb S^2$. Clearly, $X=M/D_{2m}$ is a complete Alexandrov $4$-space of curvature $\ge 0$ with a soul $S=\Bbb S^2/D_{2m}$. Note that $S$ contains three singular points of $X$: one
corresponds to the isotropy group $\Bbb Z_m$ and the other two correspond to the
isotropy group $\Bbb Z_2$.
\endexample

As shown in Example 0.5, the obstacle for Theorem A is that $F_v$ has no splitting
structure;  $F_v$ may not be locally convex (nor every direction in $\Sigma_xF_v$
points a radial curve in $F_v$, see [Pet1]), and there may not be any isometry (nor a distance
non-increasing map) from $S$ to $\tilde S_c$.

We overcome the above two obstacles in the following two theorems (see Theorem 0.6 and 0.7).

\proclaim{Theorem 0.6} Let $X$ be a complete non-compact Alexandrov space with $\op{cur}\geq0$
and an empty boundary. Assume $X$ has a soul of codimension $2$. Then

\noindent {\rm (0.6.1)} the distance function $d(\cdot, \Cal F)$ is concave on $X\setminus\Cal F$.

\noindent {\rm (0.6.2)} the gradient flow of $d(\cdot,\Cal F)$ determines a family of distance non-increasing flows from $\Cal F$, each of which flows some $F_v$ to $X\setminus\Cal F$,
provided that $\dim(X)=4$ and $X$ is simply connected and locally orientable.
\endproclaim

In verifying the concavity for $d(\cdot ,\Cal F)$, not only $F_v$ may not be locally convex,
but also the local structure of $\Cal F$ at points in $S$ is more complicated. We first establish a criterion for $d(\cdot, \Cal F)$ to be concave at $x\in X\setminus\Cal F$ in terms of properties of $(\Sigma_yX,\Sigma_y\Cal F)$ (see Lemma 2.1), where $y\in \Cal F$ satisfies $|xy|=|x\Cal F|$. If $y\notin S$, then $\Sigma_yX$ is a spherical suspension over $\Sigma_y\partial\Omega_c$, $c=f(y)$.
To verify the criterion, the key is to study the multi-valued map, $\varphi_c:  S\to \tilde S_c=F_v\cap \partial \Omega_c$ by $\varphi_c(x)=\pi^{-1}(x)\cap \tilde S_c$ (the `inverse'
of $\pi|_{\tilde S_c}: \tilde S_c\to S$ which is a branched
metric cover), whose `differential' can be defined and shares the same properties of $\varphi_c$. If $y\in S$, the verification of the criterion is more complicated, partially due to the lack of a suspension structure. Observe that $\uparrow_y^x, \Uparrow_y^{\partial \Omega_c}\, \subset  (\Sigma_yS)^\perp$ (the orthogonal complement of $\Sigma_yS$), $|\uparrow_y^x\Uparrow_y^{\partial \Omega_c}|=\frac \pi2$,
and thus $\op{dim}((\Sigma_yS)^\perp)=1$, where $\uparrow_y^x$ denotes a direction of a minimal geodesic from $y$ to $x$. Such geometric structures were studied in [RW], based on which we are able to verify the criterion case by case.

In the proof of (0.6.2), a key is to show that each $F_v\setminus S$ separates
its a neighborhood into two components.

For our purpose, without loss of generality we may assume that $X$ is also locally
orientable by [HS] (see Lemma 7.1).

\proclaim{Theorem 0.7} Let $\Psi_t|_{[0,+\infty)}:F_v\to X$ be a distance
non-increasing flow in (0.6.2).
Then for $S_c\triangleq F_v\cap \partial \Omega_c$, there is a one-parameter family of distance non-increasing onto maps, $\widetilde{\pi\circ\Psi_t}: \tilde S_c\to \tilde S_c$, with $\widetilde{\pi\circ\Psi_0}=\op{id}_{\tilde S_c}$ such that the following diagram commutes:
$$\CD \tilde S_c@>\widetilde{\pi\circ\Psi_t}>>\tilde S_c
\\ @VV \Psi_t V @ VV \pi V
\\
\Psi_t(\tilde S_c) @>\pi>>S\,\,.
\endCD$$
Consequently, $\widetilde{\pi\circ\Psi_t}$ is an isometry, and thus for all $x\in \Psi_t(\tilde S_c)$ and $r>0$,
$\pi(B_r(x))=B_r(\pi(x))$.
\endproclaim

Observe that for $c=\max f$, $\tilde S_c=S$ and $\widetilde{\pi\circ \Psi_t}=\pi\circ \Psi_t$.
For $c<\max f$, we will show that there is a finite set, $Q\subset S$,
such that $\pi: \tilde S_c\setminus\varphi_c(Q)\to S\setminus Q$ is a metric $k$-cover, and the number of $\varphi$-image, $|\varphi_c(q)|<k$ for $q\in Q$.
Moreover, $\tilde S_c$ is an Alexandrov
two sphere of non-negative curvature. In constructing $\widetilde{\pi\circ\Psi_t}: \tilde S_c\to \tilde S_c$, the geometry of $\varphi_c:
S\to \tilde S_c$ is crucial to show that for sufficiently small $t$,
$|(\pi\circ \Psi_t)^{-1}(Q)|=|\varphi_c(Q)|$; whose proof
also replies on that $S$ and $\tilde S_c$ are homeomorphic to a two sphere.
Then the following diagram,
$$\CD @. \tilde S_c\setminus\varphi_c(Q) \\ @.  @VV \pi V \\
\tilde S_c\setminus(\pi\circ \Psi_t)^{-1}(Q) @> \pi\circ \Psi_t >> S\setminus Q,
\endCD$$
satisfies that $\pi_*(\pi_1(\tilde S_c\setminus\varphi_c(Q))=(\pi\circ \Psi_t)_*(\pi_1(\tilde S_c\setminus(\pi\circ \Psi_t)^{-1}(Q)))$, which implies a lifting map of $\pi\circ \Psi_t$. Since the lifting map is distance non-increasing, it uniquely
extends to the desired map, $\widetilde{\pi\circ\Psi_t}: \tilde S_c\to \tilde S_c$, and thus $\widetilde{\pi\circ\Psi_t}$ is a distance non-increasing onto map (so is an isometry) for all $t$.
This implies that,
for all $x\in \Psi_t(\tilde S_c)$ and $r$, $\pi(B_r(x))=B_r(\pi(x))$.

Note that (A1) follows from Theorems 0.6 and 0.7. Let $x\in X\setminus\Cal F$, and let $y\in \Cal F$ such that $|xy|=|x\Cal F|$.
It turns out that the union of $x$ with $y\in \Cal F\setminus S$ is dense in $X$, and
that such a $y$ belongs to a unique $F_v$ so that the distance non-increasing flow in (0.6.2), $\Psi_t:  F_v\to X$, satisfies that $x=\Psi_{|x\Cal F|}(y)$. By Theorem 0.7, $\pi(B_r(x))=B_r(\pi(x))$, which actually holds for all $x$ by the continuity of $\pi$.

\vskip2mm

The rest of the paper is organized as follows:

\vskip2mm

In Section 1, we review basic properties of the Sharafutdinov retraction ([Per2], [Sh]), and
finite quotient of joins and radial cone-neighborhood isometries ([RW]).

From Section 2 to Section 6, we prove (0.6.1), i.e. $d(\cdot,\Cal F)$ is concave. In Section 2,
we provide a criterion for a distance function to be concave. In Section 3 and 4, the structures of $\Cal F$
and $\Sigma_y\Cal F\subset\Sigma_yX$ are described in Theorem 3.1 and 4.1 respectively.
In Section 5, we prove (3.1.3) in Theorem 3.1. And we complete the proof of (0.6.1) in Section 6.

In Section 7, we show that (A1) can be reduced to the case where $X$ is simply connected and locally orientable.

Throughout Sections  8-10, we assume that $\dim(X)=4$ and $\dim(S)=2$. In Section 8,
we show that $\tilde S_c$ is a simply connected closed
Alexandrov surface of non-negative curvature, and in Sections 9 and 10, we prove
(0.6.2) and Theorem 0.7 respectively.

\vskip4mm

We conclude the introduction with a list of symbols in the rest of the paper:

\vskip2mm

\noindent $\cdot$  $\op{Alex}^n(\kappa)$: the collection of complete Alexandrov $n$-spaces with curvature $\geq\kappa$.

\vskip1mm

\noindent $\cdot$ $B_r(p)$: a closed ball in $X$ with center $p\in X$ and radius $r$.

\vskip1mm

\noindent $\cdot$ $|A|$: the size of a subset $A\subset X$.

\vskip1mm

\noindent $\cdot$ $\partial A$: the boundary of $A\subset X$ with respect to induced topology.

\vskip1mm

\noindent $\cdot$ $A^\circ$: the union of interior points of $A$.

\vskip1mm

\noindent $\cdot$ $|xy|$: the distance between $x$ and $y$.

\vskip1mm

\noindent $\cdot$ $|xA|$: the distance from $x$ to $A$.

\vskip1mm

\noindent $\cdot$ $A^{\ge \frac \pi2}$: $\overset{\triangle}\to =\{x\in X|\, |xA|\ge \frac \pi2\}$.

\vskip1mm

\noindent $\cdot$ $A^\perp$: $\overset{\triangle}\to =A^{=\frac \pi2}\overset{\triangle}\to =\{x\in X|\, |xa|=\frac \pi2, a\in A\}$.

\vskip1mm

\noindent $\cdot$ $[xy]$: a minimal geodesic from $x$ to $y$.

\vskip1mm

\noindent $\cdot$ $[aB]$: the union of all $[ab]$ with $b\in B$.

\vskip1mm

\noindent $\cdot$ $\uparrow_x^y$: the direction of $[xy]$ at $x$.

\vskip1mm

\noindent $\cdot$ $\Sigma_aA$: the space of directions of $A$ at $a$.

\vskip1mm

\noindent $\cdot$ $\Uparrow_x^A$: the collection of directions of $[xy]$ at $x$ such that $|xy|=|xA|$ and $y\in A$.

\vskip1mm

\noindent $\cdot$ $\Bbb S^m$: the unit $m$-dimensional sphere.

\vskip1mm

\noindent $\cdot$ $A*B$: the join of $A, B\in \op{Alex}(1)$ (see Section 1).

\vskip1mm

%%%%%%%%%%%%%%%%%%%%%%%%%%%%%%%%%%%%%% section 1 Preliminaries %%%%%%%%%%%%%%%%%%%%%%%%%%

\vskip6mm

\head 1. Preliminaries
\endhead

\vskip4mm

In this section, we will supply notions and properties that will be used through the rest of the paper, and the references are [BGP], [Per2], [Pet1] and [RW].

\vskip4mm

\subhead a. Busemann functions and directions of rays
\endsubhead

\vskip4mm

In the rest of the paper, we will use $X\in \op{Alex}^n(0)$ to specify that $X$ is non-compact and with an empty boundary.
For $p\in X$, the Busemann function, $f: X\to \Bbb R$, is defined by
$$ f(x)=\lim_{t\to\infty} (|x\partial B_t(p)|-t).$$
Note that different choices of base points yield Busemann functions different by a constant (cf. [Li]). Following the approach in Riemannian case with replacing the Hessian comparison by suitable supporting functions, one gets (cf. [Per2]):

\proclaim{Lemma 1.1}  Let $X\in \op{Alex}^n(0)$. Then the Busemann function $f: X\to \Bbb R$,
 is a proper concave function.
\endproclaim

By Lemma 1.1, $f$ has a maximum $c_0$, and
any sublevel set $\Omega_c\triangleq f^{-1}([c,c_0])$ is totally convex, i.e.
%\footnote{We say that a subset $A\subseteq X$ is convex (resp. totally convex)
for any $x, y\in \Omega_c$, all minimal geodesics
jointing $x$ and $y$ lie in $\Omega_c$. Moreover, it is easy to
verify the following properties.

\proclaim{Lemma 1.2} For $p$ with $f(p)>c$, any $v\in \Uparrow_p^{\partial \Omega_c}$ is tangent to a ray at $p$, $\gamma_{p,v}$. Moreover, for all $c_1, c_2<f(p)$, $\Uparrow_p^{\partial\Omega_{c_1}}=\Uparrow_p^{\partial\Omega_{c_2}}$,
and for any $x\in\Omega_c$, $f(x)=|x\partial\Omega_c|+c$.
\endproclaim

As in the Riemannian case, using $f$ one can construct a closed totally convex subset with an
empty boundary ([Per2]), $S\subset X$, called
a soul $S$ of $X$: let $C_0=f^{-1}(c_0)$, a compact convex subset. If $\partial C_0=\emptyset$,
then $S=C_0$. Otherwise, consider $f_1(x)=|x\partial C_0|:
C_0\to \Bbb R$. Because $C_0$ is convex and $X\in \op{Alex}(0)$, $f_1$ is concave.
Let $C_1=f_1^{-1}(\max f_1)$.  If $\partial C_1=\emptyset$, then
$S=C_1$. Otherwise, repeat this process. From the decreasing sequence of integers, $\dim (X)>\dim(C_0)>\cdots$, in $k\le n$ steps one gets a totally convex set without boundary,
say $C_k=S$.

%, let $\Uparrow_p^{\partial \Omega_c}$ denote the set of directions at $p$ that tangent to a minimal geodesic $\gamma$ from $p$ to a point $\partial \Omega_c$
% i.e., the length, $|\gamma|=|p\partial \Omega_c|$.

To describe the structure of $\Sigma_pX$, $p\in S$, we need the following property.

\proclaim{Lemma 1.3 ([Ya])}  Let $Y\in \text{\rm Alex}^n(1)$, and let
$A$ be a compact (locally) convex subset in $Y$. If $\partial
A=\emptyset\ $ \footnote{Note that $A$ belongs to Alex$^m(1)$ for
some $m\leq n$, and it is our convention that $A$ consists of one
point or two points with distance $\pi$ if $m=0$. Moreover, when
$m=0$, ``$\partial A=\emptyset$'' means that $A$ consists of two
points with distance $\pi$.}, then $A^{\geq
\frac{\pi}{2}}=A^{=\frac{\pi}{2}}$.
\endproclaim

Because $S\subset X$ is convex and without boundary, $\Sigma_pS\subset \Sigma_pX$ is
convex and without boundary, where $p\in S$.
By Lemma 1.3, $(\Sigma_pS)^{\ge \frac \pi2}=(\Sigma_pS)^\perp$.
By Lemma 1.2, $\Uparrow_p^{\partial \Omega_c}\subseteq (\Sigma_pS)^{\ge \frac \pi2}$, so
$$\Uparrow_p^{\partial\Omega_{c}}\subseteq (\Sigma_pS)^\perp.\eqno{(1.1)}$$
%We point it out that without loss of generality the inclusion in (1.1) will always be
%strict; otherwise, $\pi: X\to S$ will be a submetry (see ???).
%%It is not hard to see that $\gamma(t)|_{[0,+\infty)}$ is perpendicular to $\partial \Omega_{f(\gamma(t))}$ for
%any $t>0$, and
%thus the space of directions of $X$ at $\gamma(t)$
%$$\Sigma_{\gamma(t)}X=
%\{\gamma^+(t),\gamma^-(t)\}*\Sigma_{\gamma(t)}\partial\Omega_{f(\gamma(t))}\eqno{(1.1)}$$
%(refer to Subsection 2.1 for the definition of $Y_1*Y_2$). In this
%paper, we always assume that $t$ is the arc-length parameter of a
%ray $\beta(t)|_{[0,+\infty)}$.

%%%%%%%%%%%%%%%%%%%%%%%%%%%%%%%% section 2    %%%%%%%%%%%%%%%%%%%%%
\vskip4mm

\subhead b. Sharaftdinov retraction and construction of $\Cal F$ via flat strips
\endsubhead

\vskip4mm

A distance function is semi-concave in the sense of [Pet1], and thus its gradient,
and the gradient flow, are well defined. Moreover, the gradient flow of any concave
function is distance non-increasing.

Applying to the distance functions $\{f, f_i\}$ from the construction of $S$, and patching together
the gradient flows of $\{f, f_i\}$, one obtains a one-parameter family of maps, $\phi_t: X\to X$, such that $\phi_0=\op{id}_X$ and
$\pi(x)=\lim_{t\to \infty} \phi_t(x): X\to S$, is distance non-increasing, called a Sharafutdinov retraction.

Let $$\Cal{F}\triangleq \bigcup_{p\in S,\ v\in
\Uparrow_p^{\partial\Omega_{c}}}\gamma_{v,p},$$ where
$\gamma_{v,p}$ is the ray starting from $p$ with direction $v\in
\Uparrow_p^{\partial\Omega_{c}}$. Note that any $\gamma_{v,p}$ is a gradient
curve of $f$ and thus $\pi(\gamma_{v,p})=p$.

To see that $\pi: \Cal F\to S$ is a submetry, we need the following
partial generalization of the Perel'man Flat Strip property by
Shioya-Yamaguchi ([SY])

\proclaim{Lemma 1.4} Let $X\in \op{Alex}(0)$, let $C\subset X$ be a closed convex subset
with $\partial C\neq\emptyset$. Assume a minimal geodesic $c(t)|_{[0,1]}\subset C$ satisfies
$\text{\rm dist}_{\partial C}|_{c(t)}$ is a constant, and there is a
minimal geodesic $\gamma_0$ from $c(0)$ to $\partial C$ which is
perpendicular to $c(t)$. Then there is a
minimal geodesic $\gamma_1$ from $c(1)$ to $\partial C$ such that
$\{\gamma_0, c(t), \gamma_1\}$ bounds a flat rectangle (the $4$-th
side of which belongs to $\partial C$) which is convex in $X$.
\endproclaim

%Here, that $\gamma_p$ is a minimal geodesic from $p$ to $\partial C$
%means that $\gamma_p=[px]$ with $x\in \partial C$ and
%$|px|=|p\partial C|$.

%Recall that $\Omega_c$ ($c<c_0$) is compact and totally convex in
%$X$ (see Section 1), and $S$ is totally convex in $\Omega_c$ with
%$\partial S=\emptyset$.

Lemma 1.4, together with Lemmas 1.2 and 1.3,  yields (cf. [Li]):

\proclaim{Lemma 1.5} Let $c(t)|_{[0,1]}$ be a minimal geodesic in $S$, and
let $\gamma_{v,c(0)}$ be a ray with $v\in \Uparrow_{c(0)}^{\partial \Omega_c}$.
Then there is a ray $\gamma_{w,c(1)}$ with $w\in \Uparrow_{c(1)}^{\partial\Omega_{c}}$
such that $\{\gamma_{v,c(0)},c(t), \gamma_{w,c(1)}\}$ bounds a flat strip which
is convex in $X$.
\endproclaim

By Lemma 1.5, $\pi: \Cal F\to S$  is a submetry. Moreover, Lemma 1.5 implies that
any ray in $\Cal F$ starting from $S$ coincides with some $\gamma_{v,p}$.

\remark{Remark \rm 1.6} Observe that if $\Cal F=X$, then $\pi:X\to S$ is a submetry
(so Theorem A holds). Hence, in the rest of paper we will always assume
$\Cal F\subsetneq X$. Together with (0.4.1), we may assume that $\dim(S)\le n-2$;
when $\dim(S)=n-1$, $\Cal F=X$ because $(\Sigma_pS)^\perp=\{v\}$ or $\{v_1,v_2\}$ with $|v_1v_2|=\pi$,
and thus $v$ and $ v_i$ are tangent to rays.
\endremark

\vskip4mm

\subhead c.  The rigidity of Toponogov comparison and a finite quotient of join
\endsubhead

\vskip4mm

As pointed out in Introduction 0, a verification of concavity for $d(\cdot ,\Cal F)$ requires studying
%%one obstacle in our approach to Theorem A is to verify the concavity of $d(\cdot,\Cal F)$, because $\Cal F$ may not be local convex. Given $x\in X-\Cal F$, let $q$ be the projection of $x$ on $\Cal F$.  Inspired by the construction of a supporting function ([Per]), we establish a criterion in terms of the structures, $\Sigma_q\Cal F\subsetneq \Sigma_qX$; expressing $\Sigma_qX$ in term of  $\Sigma_q\Cal F$. If $q\in \Cal F-S$, then $\Sigma_qX$ is isometric to a spherical suspension of $\Sigma_q(\Cal F\cap \partial \Omega_c)$, $c=f(q)$. The more complexity occurs in the situation $q\in S$ i.e.,
structures of $(\Sigma_pS,\Sigma_pX)$. For any $v\ne v'\in \Sigma_pS$ and $w\in (\Sigma_pS)^\perp$, any triangle $\triangle wvv'\subset \Sigma_pX$ achieves equality in the
following theorem.

\proclaim{Theorem 1.7} (Toponogove comparison) Let $Y\in\text{\rm Alex}^n(\kappa)$, and let
$\Bbb S^2_\kappa$ denote the complete simply connected $2$-manifold of curvature $\kappa$.

\noindent {\rm(1.7.1)} Given any $p\in Y$ and $[qr]\subset Y$, there is
$\tilde p$ and a $[\tilde q\tilde r]$ in $\Bbb S^2_\kappa$
with $|\tilde p\tilde q|=|pq|,|\tilde p\tilde r|=|pr|$ and $|\tilde
r\tilde q|=|rq|$ such that for any $s\in[qr]$ and $\tilde s\in[\tilde
q\tilde r]$ with $|qs|=|\tilde q\tilde s|$, $|ps|\geq|\tilde p\tilde s|$.

\vskip1mm

\noindent {\rm(1.7.2)} Given any $[qp]$ and $[qr]$ in $Y$, there is
$[\tilde q\tilde p]$ and $[\tilde q\tilde r]$ in $\Bbb S^2_\kappa$
with $|\tilde q\tilde p|=|qp|$, $|\tilde q\tilde r|=|qr|$ and
$\angle\tilde p\tilde q\tilde r=\angle pqr$ such that
$|\tilde p\tilde r|\geq|pr|$.

\vskip1mm

\noindent {\rm(1.7.3) \bf([GM2])} In {\rm (1.7.2)} (resp.
in {\rm (1.7.1)} for some interior point $s$ in $[qr]$), if equality holds, then
there exists a $[pr]$ (resp. $[qp]$ and $[pr]$) such that $[qp]$,
$[qr]$ and $[pr]$ bounds a convex surface which can be isometrically
embedded into $\Bbb S^2_\kappa$.
\endproclaim

For $p\in S$, since $(\Sigma_pS)^\perp=(\Sigma_pS)^{\geq\frac \pi2}$ (see Lemma 1.3),
$(\Sigma_pS)^\perp$ is totally convex in $\Sigma_pX$ (by (1.7.1)). A underlying geometry structure to the ``$\frac \pi2$-apart''
between $\Sigma_pS$ and $(\Sigma_pS)^\perp$, is a spherical join structure.
For convenience of readers, we recall the definition ([BGP]). For $Y_i\in
\text{Alex}^{n_i}(1)$ ($i=1, 2$, and $Y_i=\{p_i\}$ or $\{p_i,\tilde
p_i\}$ with $|p_i\tilde p_i|=\pi$ if $n_i=0$),
$$Y_1*Y_2=Y_1\times Y_2\times \left[0,\frac \pi2\right]/(y_1,y_2,0)\sim (y_1,y_2',0),
\left(y_1,y_2,\frac \pi2\right)\sim \left(y_1',y_2,\frac
\pi2\right)$$ equipped with the metric
$$\cos |(y_1,y_2,t)(y_1',y_2',t')|=\cos t\cos t'\cos |y_1y_1'|+\sin t\sin
t'\cos |y_2y_2'|$$ is called the spherical join of $Y_1$ and $Y_2$.
And $Y_1*Y_2$ is also called a spherical suspension over
$Y_2$ if $Y_1=\{p_1,\tilde p_1\}$ with $|p_1\tilde p_1|=\pi$.  It is
easy to see that $Y_1*Y_2\in \text{Alex}^n(1)$ with $n_1+n_2=n-1$
and $|p_1p_2|=\frac \pi2$ for all $p_i\in Y_i$, and there is a
unique minimal geodesic jointing $p_1$ with $p_2$.

Observe that if the following property hods that for $v\in
\Sigma_pS$ and $w\in (\Sigma_pS)^\perp$, there is a unique minimal geodesic from
$v$ to $w$, then there is an isometric embedding, $\Sigma_pS*(\Sigma_pS)^\perp\hookrightarrow \Sigma_pX$ ([Li]).
Unfortunately, if $X$ is not a topological manifold around $S$, then one cannot expect such
join structure in $\Sigma_pX$.

In [RW], we classify $X\in \op{Alex}^n(1)$ which contains two convex subsets
without boundary such that the sum of dimension is $n-1$. Applying this result to
the case $\dim(S)=n-2$ and $(\Sigma_pS)^\perp$ is a circle, we conclude $\Sigma_pX$ as a finite
quotient of join (see (1.8.2) below).

\proclaim{Theorem 1.8 ([RW])} Let $Y\in \text{\rm Alex}^n(1)$, and
let $Y_1, Y_2$ be two compact convex subsets in $Y$ such that
$|y_1y_2|=\frac \pi2$ for any $y_i\in Y_i$. Then the following
holds.

\noindent {\rm (1.8.1)} $n_1+n_2\le n-1$, where $n_i=\dim(Y_i)$.

\noindent {\rm (1.8.2)}  If $n_1+n_2=n-1$ and if $\partial
Y_1=\partial Y_2=\emptyset$, then $Y$ is isometric to a finite
quotient of join. In detail, there is $p_i\in Y_i$ ($i=1, 2$) and a
finite group $\Gamma$ which acts effectively and isometrically on
$(\Sigma_{p_i}Y_i)^\perp$ such that $Y_1\overset{\text{iso}}\to
\cong(\Sigma_{p_2}Y_2)^\perp/\Gamma,\ Y_2 \overset{\text{iso}}\to
\cong(\Sigma_{p_1}Y_1)^\perp/\Gamma$ and
$Y\overset{\text{iso}}\to\cong [(\Sigma_{p_1}Y_1)^\perp *
(\Sigma_{p_2}Y_2)^\perp] /\Gamma\ $ \footnote{Note that the $\Gamma$-actions on
$(\Sigma_{p_i}Y_i)^\perp$ can be extended uniquely to a
$\Gamma$-action on
$(\Sigma_{p_1}Y_1)^\perp*(\Sigma_{p_2}Y_2)^\perp$.}.
\endproclaim

%\remark{Remark \rm 1.9} In fact, if $\dim(S)=n-1$ , then we can
%prove that $X=\Cal{F}$, so $\pi$ is a submetry ([SY]). Hence, in
%studying the Soul Conjecture in Axexandrov geometry we need only to
%consider the case where
%$$\dim(S)\leq\dim(X)-2 \text{ and } X\neq\Cal F.$$
%\endremark

\vskip4mm

\subhead d.  The radial cone-neighborhood isometry
\endsubhead

\vskip4mm

The notion of a radial cone-neighborhood isometry refers to a multi-valued map, $f: Y\to \tilde Y$,
that is distance non-decreasing and preserves the metric in the radial direction. A simple example is the inverse map of the projection:
$P: M\to M/\Gamma$, where $M$ is a complete Alexadrov space and $\Gamma$ is a finite group of
isometries.

\example{Definition 1.9 ([RW])} Let $Y, \, \tilde Y\in
\text{Alex}(\kappa)$ with $\dim(\tilde Y)\ge\dim(Y)\geq1$. A
multi-valued map $f:Y\rightarrow \tilde Y$ is called a {\it radial
cone-neighborhood isometry} if the following properties hold:

\noindent(1.9.1) $m=\max_{p\in Y}\{|f(p)|\}<\infty$ and $Y_m\triangleq\{p\in Y|\
|f(p)|=m\}$ is dense in $Y$;

\noindent(1.9.2) For $p, q\in Y$, $\tilde p\in f(p)$,
$\tilde q\in f(q)$, $|\tilde p\tilde q|\geq|pq|$. Moreover,
given $[pq]$ and $\tilde p$, $\tilde q$ and $[\tilde p\tilde q]$ can
be chosen so that  $f:[pq]\rightarrow [\tilde p\tilde q]\subseteq f([pq])$ is an isometry.
\endexample

Note that (1.9.1) and (1.9.2) imply that ([RW]):

\noindent(1.9.3) For $p\in Y_m$ and $\tilde p\in f(p)$, there is
an $\epsilon_p>0$ such that
$f:B_{\epsilon_p}(p)\rightarrow B_{\epsilon_p}(\tilde p)$ is an isometric embedding.

\noindent (1.9.4) For $p\in  Y\setminus Y_m$ and $\tilde p\in f(p)$, there is
an $\epsilon_p>0$ such that $B_{\epsilon_p}(\tilde p)\cap f(Y)\subset f(B_{\epsilon_p}(p))$
and for any $[pq]\subset B_{\epsilon_p}(p)$, $f([pq])=\bigcup_{i=1}^l[\tilde p\tilde q_i]$ and
$f:[pq]\rightarrow [\tilde p\tilde q_i]$ is an isometry; moreover, if $q\in Y_m$, then
there is $\delta>0$ such that $V_{[pq],\delta}\setminus\{p\}\subset Y_m$
and that $f:V_{[pq],\delta}\rightarrow V_{[\tilde p\tilde q_i],\delta}$ is an isometric embedding, where
$$V_{[pq],\delta}\triangleq\{q'\in S
|\ |pq'|\leq |pq|,\ \exists\ [pq'] \text{ s.t. }
|\uparrow_p^{q}\uparrow_p^{q'}|<\delta\}.$$

\vskip2mm

A radial cone-neighborhood isometry has a well defined `differential'. For
$p\in Y$ and $\tilde p\in f(p)$, by (1.9.2) one can
define a multi-valued map, $\text{D}f: \Sigma'_pY\to \Sigma_{\tilde p}\tilde Y$,
$\text{D}f(\uparrow_p^q)=\{\uparrow_{\tilde p}^{\tilde q}|\, \, \tilde q\in f(q)\}$,
where $\Sigma_p'Y$ denotes the set of directions of non-trivial geodesics from
$p$ (and $\Sigma_pY$ is the closure of $\Sigma_p'Y$).

\proclaim{Proposition 1.10 ([RW])} Let $f:Y\rightarrow \tilde Y$ be a radial
cone-neighborhood isometry. Then

\noindent{\rm (1.10.1)} $\text{\rm D}f$ extends to a multi-valued map, $\text{\rm D}f:\Sigma_pY\rightarrow
\Sigma_{\tilde p}\tilde Y$, $\tilde p\in f(p)$, which is again
a radial cone-neighborhood isometry.

\noindent{\rm (1.10.2)} If $Y$ is compact with
$\partial Y=\emptyset$ and $\dim(Y)=\dim(\tilde Y)$,
then $f$ is surjective (thus $\tilde Y$ is compact) and $\partial \tilde Y=\emptyset$.
\endproclaim

We point it out that in the proof of Theorem 1.8, fixing $p\in Y_2$, we construct
a natural radial cone-neighborhood isometry,  $\phi_p:  Y_1\to (\Sigma_pY_2)^\perp$,
which is crucial in
the proof. In the present paper, there is a natural multi-valued map, $\varphi_c: S\to X$
with $\varphi_c(S)\subset \partial \Omega_c$, which is a radial cone-neighborhood isometry (see Section 4). Indeed, techniques
developed in [RW] in analyzing a radial cone-neighborhood isometry are basic tools
in our proofs of Theorems 0.6 and 0.7.

%%%%%%%%%%%%%%%%%%%%%%%%%%%%%%%% section 2   %%%%%%%%%%%%%%%%%%%%%

\vskip4mm

\head 2. The Concavity of $d(\cdot,\Cal F)$
\endhead

\vskip4mm

Our main effort in Sections 2-6 is to prove (0.6.1), where main references are [BGP] and [RW].

Given $Y\in \op{Alex}(0)$, a closed subset $Z\subset Y$, and $x\in Y\setminus Z$, let $\alpha(t)|_{[0,\epsilon)}$ be an
arc-length parameter minimal geodesic in $Y\setminus Z$ with
$\alpha(0)=x$, and let $y\in Z$ such that $|xy|=|xZ|$.
Then there is a minimal geodesic $[xy]$ such that
$$\frac{\text{\rm d}|y\alpha(t)|}{\text{\rm d}t}|_{t=0^+}=-\cos|
\uparrow_x^y\alpha^+(0)|;$$
and there is $[\alpha(t)y]$ which converges to $[xy]$ as $t\to 0$.
Note that $d(\cdot, Z)$ is concave on $Y\setminus Z$ if
for any $x\in Y\setminus Z$ and such an $\alpha(t)$
$$|\alpha(t)Z|\leq|xy|-t\cos|\uparrow_x^y\alpha^+(0)|+o(t^2).$$
For our purpose, we need the following criterion for $d(\cdot, Z)$ to be a concave function.

\proclaim{Lemma 2.1} Let $Y, Z, x, y, \alpha(t), [xy]$ and $[\alpha(t)y]$ be as in
the above. Suppose that for any small $t$, there is $w\in \Sigma_yZ\ $\footnote{Here, $\Sigma_yZ$ is defined to
be the union of all limits of $\uparrow_y^{z}$ with $z\in Z$
converging to $y$.} such that

\noindent$(2.1.1)$ $|\uparrow_y^x w|=\frac\pi2$ and $\uparrow_y^{\alpha(t)}$ lies in a $[\uparrow_y^x
w]\subset\Sigma_yY$, and either

\noindent $(2.1.2)$ $w$ points a radial curve in $Z$, or

\noindent $(2.1.3)$ there is a $w'\in\Sigma_yZ$ satisfying the following conditions:

\noindent $(2.1.3.1)$ $|\uparrow_y^x w'|=\frac\pi2$ and $|ww'|\to 0$ as $t\to 0$,

\noindent $(2.1.3.2)$ $[\uparrow_y^x w]$ and $w'$ determines a convex
spherical surface (in $\Sigma_yY$),

\noindent $(2.1.3.3)$ there is $[yz]\subset Z$ such that $|yz|\geq
t$ and $\uparrow_y^z=w'$.

\noindent Then we have that
$|\alpha(t)Z|\leq|xy|-t\cos|\uparrow_x^y\alpha^+(0)|+o(t^2).$
\endproclaim

Note that (2.1.1) and (2.1.2) are known results
([Per2], [Pet1]).

We point out that in our circumstances, $Z=\Cal F$ may not satisfy
(2.1.2); e.g., $\Cal F$ may not be locally
convex. A verification of (2.1.3)
for $\Cal F$ is based on structures of $(\Sigma_yF_v,\Sigma_yX)$ (see  Section 4), which
in turn, relies on the structure of $\Cal F$ (see Section 3).

\demo{Proof of Lemma 2.1}

As seen before the proof, it suffices to prove Lemma 2.1 with (2.1.1) and
(2.1.3). By (2.1.1) and (2.1.3.1-2) (an embedded spherical triangle
$\triangle\uparrow_y^xww'$, we get
$$\cos|\uparrow_y^{\alpha(t)}w'|=
\cos|\uparrow_y^{\alpha(t)}w|\cdot\cos|ww'|,\eqno{(2.1)}$$
and
$$|\uparrow_y^{\alpha(t)}w'|<
\frac{\pi}{2}-(1-\chi(t))\cdot|\uparrow_y^{\alpha(t)}\uparrow_y^{x}|,\eqno{(2.2)}$$
where $\chi(t)\to 0$ as $t\to0$.

Let $\triangle \tilde x\tilde y\tilde\alpha(t)\subset \Bbb R^2$ denote a
comparison triangle of $\triangle xy\alpha(t)$, and let $\tilde y'\in \Bbb R^2$
such that $[\tilde y\tilde y']$ is
perpendicular to $[\tilde y\tilde x]$ and $[\tilde y'\tilde
\alpha(t)]$ is parallel to $[\tilde y\tilde x]$. Note that $|\tilde
y\tilde y'|\leq t$, so there is  $y'\in [yz]$ such that
$|yy'|=|\tilde y\tilde y'|$ (see (2.1.3.3)). Moreover, by (1.7.2) we
have that
$$\angle\tilde\alpha(t)\tilde y\tilde y'\geq
|\uparrow_y^{\alpha(t)} w|.\eqno{(2.3)}$$
Note that if $\angle\tilde\alpha(t)\tilde y\tilde y'\geq
|\uparrow_y^{\alpha(t)} w'|,$ then Lemma 2.1 holds because
by (1.7.2) it is not hard to see
that
$$|\alpha(t)y'|\leq|\tilde\alpha(t)\tilde y'|=
|xy|-t\cos|\uparrow_{\tilde x}^{\tilde y}\uparrow_{\tilde x}^{\tilde
\alpha(t)}|\leq |xy|-t\cos|\uparrow_x^y\alpha^+(0)|.\eqno{(2.4)}$$
Observe that $\angle\tilde\alpha(t)\tilde y\tilde y'\geq
|\uparrow_y^{\alpha(t)} w'|$ is satisfied when $|\uparrow_y^x\uparrow_y^{\alpha(t)}|>2|\uparrow_{\tilde
y}^{\tilde x}\uparrow_{\tilde y}^{\tilde \alpha(t)}|$, because
by (2.2) we conclude that
$\angle\tilde\alpha(t)\tilde y\tilde y'=\frac\pi2-
|\uparrow_{\tilde y}^{\tilde x}\uparrow_{\tilde y}^{\tilde
\alpha(t)}|>|\uparrow_y^{\alpha(t)} w'| \ \text{ (as $t\to 0$)}$.

The remaining case is where
$|\uparrow_y^x\uparrow_y^{\alpha(t)}|\leq 2|\uparrow_{\tilde
y}^{\tilde x}\uparrow_{\tilde y}^{\tilde \alpha(t)}|$ (and
$\angle\tilde\alpha(t)\tilde y\tilde y'<|\uparrow_y^{\alpha(t)}
w'|$). Note that there is a constant $C$ depending on $|xy|$ and
$|\uparrow_x^y\alpha^+(0)|$ such that $|\uparrow_{\tilde y}^{\tilde
x}\uparrow_{\tilde y}^{\tilde \alpha(t)}|<C\cdot t$, so
$$|\uparrow_{y}^{x}\uparrow_{y}^{\alpha(t)}|<2C\cdot t.\eqno{(2.5)}$$
In this case, we can select $\tilde y''$ in the plane containing
$\triangle \tilde x\tilde y\tilde\alpha(t)$ such that $\angle
\tilde\alpha(t)\tilde y\tilde y''=|\uparrow_y^{\alpha(t)} w'|$ and
$[\tilde y\tilde y'']$ is perpendicular to $[\tilde y''\tilde
\alpha(t)]$. Note that $|\tilde y\tilde y''|<|\tilde y\tilde y'|$,
so we can find $y''\in [yz]$ such that $|yy''|=|\tilde y\tilde
y''|$. Then we have that
$$\align |\alpha(t)y''|&\overset{\text{by (1.7.2)}}\to\leq|\tilde\alpha(t)\tilde y''|\\
&=|\tilde\alpha(t)\tilde y'|+|\tilde\alpha(t)\tilde y''|-|\tilde\alpha(t)\tilde y'|\\
&=|\tilde\alpha(t)\tilde y'|+|\tilde
y\tilde\alpha(t)|\left(\sqrt{1-\cos^2|\uparrow_y^{\alpha(t)}w'|}-
\sqrt{1-\cos^2|\uparrow_{\tilde y}^{\tilde \alpha(t)}\uparrow_{\tilde y}^{\tilde y'}|}\right)\\
&=|\tilde\alpha(t)\tilde y'|+|\tilde y\tilde\alpha(t)|
\frac{\cos^2|\uparrow_{\tilde y}^{\tilde \alpha(t)}\uparrow_{\tilde
y}^{\tilde y'}|-\cos^2|\uparrow_y^{\alpha(t)}w'|}
{\sqrt{1-\cos^2|\uparrow_y^{\alpha(t)}w'|}+
\sqrt{1-\cos^2|\uparrow_{\tilde y}^{\tilde
\alpha(t)}\uparrow_{\tilde y}^{\tilde y'}|}}\\
&\overset{\text{by
(2.1)}}\to=|\tilde\alpha(t)\tilde y'|+|\tilde y\tilde\alpha(t)|
\frac{\cos^2|\uparrow_{\tilde y}^{\tilde \alpha(t)}\uparrow_{\tilde
y}^{\tilde y'}|-\cos^2|\uparrow_y^{\alpha(t)}w|\cdot\cos^2|ww'|}
{\sqrt{1-\cos^2|\uparrow_y^{\alpha(t)}w'|}+
\sqrt{1-\cos^2|\uparrow_{\tilde y}^{\tilde
\alpha(t)}\uparrow_{\tilde y}^{\tilde y'}|}}\\
\endalign$$
$$\align
&\overset{\text{by (2.3)}}\to\leq|\tilde\alpha(t)\tilde y'|+|\tilde
y\tilde\alpha(t)|
\frac{\cos^2|\uparrow_y^{\alpha(t)}w|(1-\cos^2|ww'|)}
{\sqrt{1-\cos^2|\uparrow_y^{\alpha(t)}w'|}+
\sqrt{1-\cos^2|\uparrow_{\tilde y}^{\tilde
\alpha(t)}\uparrow_{\tilde y}^{\tilde y'}|}}\\
&=|\tilde\alpha(t)\tilde y'|+|\tilde y\tilde\alpha(t)|
\frac{\sin^2|\uparrow_y^x\uparrow_y^{\alpha(t)}|(1-\cos^2|ww'|)}
{\sqrt{1-\cos^2|\uparrow_y^{\alpha(t)}w'|}+
\sqrt{1-\cos^2|\uparrow_{\tilde y}^{\tilde
\alpha(t)}\uparrow_{\tilde y}^{\tilde y'}|}}\\
&\overset{\text{by (2.5)}}\to\leq|\tilde\alpha(t)\tilde y'|+o(t^2).
\endalign$$
We thereby complete the proof by taking (2.4) into account.
\hfill$\qed$
\enddemo

%%%%%%%%%%%%%%%%%%%%%%%%%%%%%%%% section 4 %%%%%%%%%%%%%%%%%%%%%

\vskip4mm

\head 3. Structures of $\Cal F\subsetneq X$ with $\dim(S)=\dim(X)-2$
\endhead

\vskip4mm

By Lemma 2.1, $d(\cdot, \Cal F)$ is concave if $\Cal F$ satisfies (2.1.1) and (2.1.3).
Our verification of (2.1.1) and (2.1.3) is divided into
three steps: we first analyze structures of $\Cal F$ in Sections 3 and 4, based on which we then analyze
structures on $(\Sigma_y\Cal F,\Sigma_yX)$ in Section 4.
In Section 6, we show that $\Cal F$ satisfies (2.1.1) and (2.1.3).
Note that in Sections 3-6, $\dim(\Cal F)\ge \dim(X)-1$, $\dim(S)=\dim(X)-2$  and $\dim(X)\ge 4$.

Let $S_0\subset S$ consisting of points whose spaces of directions are isometric to
a unit sphere. Then $S_0$ is totally convex in $S\in \op{Alex}(0)$,
with a full measure ([OS]). Fixing $p_0\in S_0$ and $v\in \Uparrow_{p_0}^{\partial
\Omega_c}$, we define a subset, $\overset{\circ}\to F_v\subset\Cal F$, as follows.
For any $p\in S_0$ and $[p_0p]$, by Lemma 1.5 there is a $w\in \Uparrow_p^{\partial \Omega_c}$ such that
$\{\gamma_{v,p_0},[p_0p],\gamma_{w,p}\}$ bounds a flat strip. We will call $\gamma_{w,p}$
parallel to $\gamma_{p_0,v}$ along $[p_0p]$. Let
$$\split \overset{\circ}\to F_v&=\{\gamma_{w,p}|\ p\in S_0, w\in \Uparrow_p^{\partial \Omega_c}, \gamma_{p,w} \text{ is parallel to
$\gamma_{p_0,v}$ along a piecewise}\\ & \hskip16mm \text{ minimal geodesic in $S$ from $p_0$ to $p$}\}.\endsplit$$
Let $F_v$ denote the closure of $\overset{\circ}\to F_v$ in $X$. Because $\overset{\circ}\to F_v$ consists of rays, so does $F_v$.
Given $p\in S_0$, $q\in S$ and $[qp]$, and $w\in \Uparrow_q^{\partial \Omega_c}$, similar to
the above by Lemma 1.5, there is $v\in \Uparrow_p^{\partial \Omega_c}$ such that
$\gamma_{v,p}$ is parallel to $\gamma_{w,q}$ along $[qp]$. Plus the fact that $S_0$ is dense in $S$, we see

$$ \Cal F =\bigcup_{w\in \Uparrow_{p_0}^{\partial\Omega_c}}F_w.$$

We now are ready to describe basic structures on $F_v$.

\proclaim{Theorem 3.1}  Let $X\in \op{Alex}(0)$ with a soul $S$ of codimension $2$,
and let $S_0\subset S$ and $F_v\subset \Cal F$ be as in the above. Then

\noindent {\rm (3.1.1)}  Given $p, q\in S_0$, $[pq]$ and $v\in \Uparrow_p^{\partial \Omega_c}$, there is a unique $w\in \Uparrow_q^{\partial \Omega_c}$ such that $\{\gamma_{v,p},[pq],\gamma_{w,q}\}$ bounds a flat strip. Moreover, the induced bijection, $\phi_{[pq]}: \, \Uparrow_p^{\partial \Omega_c}\to \Uparrow_q^{\partial \Omega_c}$
by $v\mapsto w$, is an isometry.

\noindent {\rm (3.1.2)}  There is $k<\infty$ such that for all $p\in S_0$, $q\in S$,
$|\Uparrow_q^{\partial \Omega_c}\cap \Sigma_qF_v|\le
|\Uparrow_p^{\partial \Omega_c}\cap \Sigma_pF_v|=k$.

\noindent {\rm (3.1.3)} If $F_v\neq F_w$, then $F_v\cap F_w=S$.
\endproclaim

(3.1.1) easily implies that restricted to $\Sigma_pF_v$, $\phi_{[pq]}$ is an isometry to $\Sigma_qF_v$; (3.1.2) says that the number of rays in $F_v$ at $p\in S_0$ equals to $k$ which is independent of $p$, and the number of rays in $F_v$ at $q\in S\setminus S_0$ is bounded above by $k$; and (3.1.3) says that either $F_v=F_w$ or $F_v\cap F_w=S$.

We point out that (3.1.1) holds for $S$ of any codimension, and because our proof of (3.1.3)
relies on structures of $\Sigma_pF_v$, $p\in S$, which will be
studied in Section 4, we postpone a proof of (3.1.3) in Section 5.

We now give a proof of (3.1.1). Note that for $p\in S_0$, we know that (cf. [GW])
$$\Sigma_pX=(\Sigma_pS)*Y=\Bbb S^{m}* Y, \eqno{(3.1)}$$
where $m=\dim(S)-1$ and $Y$
belongs to Alex$(1)$ of dimension $\dim(X)-2-m$. Note that $\partial \Sigma_pX=\emptyset$
because $\partial X=\emptyset$, which implies $\partial
Y=\emptyset$, so $Y$ is a circle $S^1$ with perimeter
$\leq 2\pi$ if $\dim(S)=\dim(X)-2$.

\demo {Proof of (3.1.1)}

By Lemma 1.5, there is $w\in \Uparrow_q^{\partial \Omega_c}$ such that
$\{\gamma_{v,p},[pq],\gamma_{w,q}\}$ bounds a flat strip. Note that this flat
strip determines a minimal geodesic of length $\frac\pi2$ from $v$ to $\uparrow_p^q$ in $\Sigma_pX$.
By (3.1), the minimal geodesic from $v$ to $\uparrow_p^q$ is unique,
so is $w$ unique. Then we can define a natural map
$\phi_{[pq]}:\ \Uparrow_p^{\partial\Omega_{c}}\to\Uparrow_q^{\partial\Omega_{c}}$
with $\phi_{[pq]}(v)=w$. By taking $\phi_{[qp]}\ (=\phi_{[pq]}^{-1})$ into account,
we see that $\phi_{[pq]}$ is a bijection.

Next we will prove that $\phi_{[pq]}$ is an isometry.
(The proof is an alternative one of the proof of Lemma
2.2 in [Li].)
Let $u,v\in \Uparrow_p^{\partial\Omega_{c}}$, and let
$\bar u=\phi_{[pq]}(u), \bar v=\phi_{[pq]}(v)$.
It suffices to show that
$$|\bar u\bar v|\leq|uv|. \eqno{(3.2)}$$
Let $\gamma_u(t)|_{[0,+\infty)}$ be the (arc-length parameter)
ray with $\gamma_u(0)=p$ and $\gamma_u^+(0)=u$. By the definition of
$\phi_{[pq]}$, $\{\gamma_u(t), [pq], \gamma_{\bar u}(t)\}$ and
$\{\gamma_v(t), [pq], \gamma_{\bar v}(t)\}$ bound two flat strips
$\Cal S_u$ and $\Cal S_v$ respectively. To those, we associate two
flat strips $\tilde \Cal S_u$ and $\tilde \Cal S_v$ in Euclidean
space $\Bbb R^3$ bounded by $\{\tilde \gamma_u(t), [\tilde p\tilde
q], \tilde \gamma_{\bar u}(t)\}$ and $\{\tilde \gamma_v(t), [\tilde
p\tilde q], \tilde \gamma_{\bar v}(t)\}$ respectively, where
$|\tilde p\tilde q|=|pq|$, and the four $\tilde\gamma$ are all
perpendicular to $[\tilde p\tilde q]$, and
$$|\tilde \gamma^+_u(0)\tilde \gamma^+_v(0)|=|\tilde \gamma^+_{\bar u}(0)\tilde \gamma^+_{\bar
v}(0)|=|\bar u\bar v|.\eqno{(3.3)}$$ Hence, in order to see (3.2) it suffices to show that
$$|uv|\geq |\tilde \gamma^+_{u}(0)\tilde \gamma^+_{v}(0)|.\eqno{(3.4)}$$
By (3.1), $\Cal S_u$ and $\Cal S_v$ (similarly for $\tilde \Cal S_u$
and $\tilde \Cal S_v$) determines a geodesic triangle $\triangle
\uparrow_p^quv$  in $\Sigma_pX$ with
$|\uparrow_p^qu|=|\uparrow_p^qv|=\frac\pi2$, which can be embedded
isometrically into $\Bbb S^2$. Note that
$\uparrow_p^{\gamma_{\bar u}(t)}\in [\uparrow_p^qu]$,
$\uparrow_p^{\gamma_{\bar v}(t)}\in [\uparrow_p^qv]$,
$\uparrow_{\tilde p}^{\tilde \gamma_{\bar u}(t)}\in
[\uparrow_{\tilde p}^{\tilde q}\tilde \gamma_{u}^+(0)]$ and
$\uparrow_{\tilde p}^{\tilde \gamma_{\bar v}(t)}\in
[\uparrow_{\tilde p}^{\tilde q}\tilde \gamma_{v}^+(0)]$, and
$$|\uparrow_p^{\gamma_{\bar u}(t)}\uparrow_p^q|=|\uparrow_p^{\gamma_{\bar v}(t)}\uparrow_p^q|=
|\uparrow_{\tilde p}^{\tilde \gamma_{\bar u}(t)}\uparrow_{\tilde
p}^{\tilde q}|= |\uparrow_{\tilde p}^{\tilde\gamma_{\bar
v}(t)}\uparrow_{\tilde p}^{\tilde q}|=O(t) \text{ as } t\to
0.\eqno{(3.5)}$$ On the other hand, (3.3) implies that (ref. [BGP])
$$|\tilde\gamma_{\bar u}(t)\tilde\gamma_{\bar v}(t)|
=|\gamma_{\bar u}(t)\gamma_{\bar v}(t)|+o(t) \text{ as } t\to 0,$$
which together with (1.7.2) implies that
$$|\uparrow_p^{\gamma_{\bar u}(t)}\uparrow_p^{\gamma_{\bar v}(t)}|\geq
|\uparrow_{\tilde p}^{\tilde \gamma_{\bar u}(t)}\uparrow_{\tilde
p}^{\tilde \gamma_{\bar v}(t)}|+o(t) \text{ as } t\to
0.\eqno{(3.6)}$$ It is not hard to see that (3.5) and (3.6)
together implies (3.4). \hfill$\qed$
\enddemo

\remark{Remark \rm 3.2} We point it out that in the above proof (e.g., the proof for (3.4)),
the key ingredient is the joint structure on $\Sigma_pX\ (p\in S_0)$
(see (3.1)). Similarly,
for any $[pq]\subset S$ with $p\in S_0$, we can define a map
$$\phi_{[pq]}:\ \Uparrow_p^{\partial\Omega_{c}}\to\Uparrow_q^{\partial\Omega_{c}}
\text{ which is distance non-increasing}. \eqno{(3.7)}$$
\endremark

\proclaim{Corollary 3.3} ``$\Cal{F}\subsetneq X$'' is equivalent to that
$\Uparrow_{p}^{\partial\Omega_{c}}\subsetneq(\Sigma_{p}S)^\perp$ for any
fixed $p\in S_0$.
\endproclaim

Corollary 3.3 amounts to that if
$\Uparrow_{p}^{\partial\Omega_{c}}=(\Sigma_{p}S)^\perp$ for some
$p\in S_0$, then $X=\Cal F$ (and thus $\pi$ is a submetry, see Remark 1.6).

\demo{Proof} Assume that
$\Uparrow_{\bar p}^{\partial\Omega_{c}}=(\Sigma_{\bar p}S)^\perp$ for some
$\bar p\in S_0$. We need to show that $X=\Cal F$. By (3.1.1) and (3.1), we first
conclude that $\Uparrow_p^{\partial\Omega_{c}}=\Uparrow_{\bar p}^{\partial\Omega_{c}}$ and $\Uparrow_p^{\partial\Omega_{c}}=(\Sigma_pS)^\perp$ for
all $p\in S_0$. It follows that $$X\setminus\Cal F\subset \Cal
G\triangleq \{x\in X|\ \exists\ q\in S\setminus S_0 \text{ s.t. }
|xq|=|xS|\}.$$ Note that
$\dim((\Sigma_qS)^\perp)\leq\dim(X)-1-\dim(S)$ for any $q\in S$ (by
(1.8.1)). It then is not hard to see that $\Cal G$ has measure 0 \
\ \footnote{An essential reason here is that the distance function
to $S$ is semi-concave on $X\setminus S$ ([Pet1]). And in this
paper, ``measure'' always means the Hausdorff measure.} because
$S\setminus S_0$ has measure 0 on $S$. On the other hand, it is easy
to see that $\Cal F$ is closed in $X$ (by Lemma 1.5), so $X\setminus
\Cal F$ must have a positive measure if $\Cal{F}\subsetneq X$. Hence, it has
to hold that $X=\Cal F$. \hfill$\qed$
\enddemo

\demo{Proof of (3.1.2)}

By (3.1.1), it is easy to see
that, for all $p\in S_0$, $\Uparrow_p^{\partial \Omega_c}\cap \Sigma_p
\overset{\circ}\to F_v$ is isometric to $\Uparrow_{p_0}^{\partial \Omega_c}
\cap \Sigma_{p_0}\overset{\circ}\to F_v$. Because $F_v$ is closure of
$\overset{\circ}\to F_v$, for $q\in S\setminus S_0$, $|\Uparrow_q^{\partial\Omega_c}\cap \Sigma_qF_v|\le
|\Uparrow_{p_0}^{\partial \Omega_c}\cap \Sigma_{p_0}\overset{\circ}\to F_v|$. Hence, it remains to show
that $|\Uparrow_{p_0}^{\partial \Omega_c}\cap \Sigma_{p_0}\overset{\circ}
\to F_v|<\infty$.

Let $v\neq \bar v\in \Uparrow_{p_0}^{\partial \Omega_c}\cap \Sigma_{p_0}\overset{\circ}\to F_v$.
By the definition of $\overset{\circ}\to F_v$,
$\gamma_{\bar v,p_0}$ is parallel to $\gamma_{v,p_0}$ along a
piecewise minimal geodesic
$\bigcup_{j=1}^l[p_jp_{j+1}]\subset S_0$ with $p_1=p_{l+1}=p_0$. By
(3.1.1), the isometry
$$\phi_{[p_lp_{l+1}]}\circ\phi_{[p_{l-1}p_l]}\circ\cdots\circ\phi_{[p_1p_2]}:\
\Uparrow_{p_0}^{\partial\Omega_{c}}\to\Uparrow_{p_0}^{\partial\Omega_{c}}\eqno{(3.8)}$$
maps $v$ to $\bar v$. Since we have assumed that $S$ has
codimension 2,
$$(\Sigma_{p_0}S)^\perp \text{ is a
circle (see (3.1))}.$$ Then the isometry in (3.8) must
be an isometry of $(\Sigma_{p_0}S)^\perp$ restricted to
$\Uparrow_{p_0}^{\partial\Omega_{c}}$ (cf. the proof of Key Lemma 0.7
in [Li]). Thereby, if
$|\Uparrow_{p_0}^{\partial \Omega_c}\cap \Sigma_{p_0}\overset{\circ}\to F_v|=\infty$, we can conclude
that $\Uparrow_{p_0}^{\partial\Omega_{c}}=(\Sigma_{p_0}S)^\perp$, which
contradicts Corollary 3.3 because we have assumed that $\Cal{F}\subsetneq X$.
\hfill$\qed$
\enddemo

\remark{Remark \rm 3.4} For $q\in S$, let
$\xi_1,\xi_2\in \Uparrow_{q}^{\partial\Omega_{c}}$ such that
$\gamma_{\xi_i,q}\subset F_v$. And for small ball $B_\epsilon(q)\subset S$,
let $F_v|_{B_\epsilon(q),\xi_i}$
denote the collection of rays in $F_v$ each of which is parallel to
$\gamma_{\xi_i,q}$ along a minimal geodesic in
$B_\epsilon(q)$ starting from $q$. Note that (3.1.2) implies
that for small $\epsilon>0$,
$$F_v|_{B_\epsilon(q),\xi_1}\cap
F_v|_{B_\epsilon(q),\xi_2}=B_\epsilon(q)\eqno{(3.9)}$$  and
if $\gamma_{w,q'}\subset F_v$ with $q'\in B_\epsilon(q)$
is parallel to a ray in $F_v|_{B_\epsilon(q),\xi_i}$, then
$$\gamma_{w,q'}\subset F_v|_{B_\epsilon(q),\xi_i}.\eqno{(3.10)}$$
By (3.9) it is easy to see that
$$ S_k(F_v)\triangleq\{q\in S\,|\, |\Uparrow_{q}^{\partial \Omega_c}\cap \Sigma_{q} F_v|=k,\
\text{i.e. at $q$ there are
$k$ pieces of rays in $F_v$}\}$$ is open and dense in $S$.
Note that
$S_0\subseteq S_k(F_v)$, and that for any
$[qp]\subset S$ with $p\in  S_k(F_v)$,
$[qp]\setminus\{q\}\subset S_k(F_v)$.
And (3.9) and (3.10) imply that
$$\text{$\gamma_{w,q'}\subset F_v|_{B_\epsilon(q),\xi_i}$ is parallel to $\gamma_{\xi_i,q}$ along any $[qq']$},\eqno{(3.11)}$$
and thus there is $l\leq k$ such that for all $q'\in B_\epsilon(q)$ and $q''\in B_\epsilon(q)\cap  S_k(F_v)$,
$$\left|\Uparrow_{q'}^{\partial \Omega_c}\cap \Sigma_{q'} F_v|_{B_\epsilon(q),\xi_i}\right|\leq
\left|\Uparrow_{q''}^{\partial \Omega_c}\cap \Sigma_{q''} F_v|_{B_\epsilon(q),\xi_i}\right|=l.\eqno{(3.12)}$$
\endremark

%%%%%%%%%%%%%%%%%%%%%%%%%%%%%%%% section 5 %%%%%%%%%%%%%%%%%%%%%

\vskip4mm

\head 4. Radial Cone-Neighborhood Isometries and Structures on $(\Sigma_yF_v, \Sigma_yX)$.
\endhead

\vskip4mm

This section is a part of preparation for Section 6, where we show
that $d(\cdot,\Cal F)$ satisfies (2.1.1) and (2.1.3). We will analyze structures
on $(\Sigma_y\Cal F,\Sigma_yX)$; and if $y\in \Cal F\setminus S$,
it reduces to analyze structures on $(\Sigma_yF_v,\Sigma_yX)$, because
by (3.1.3) there is a unique $F_v\ni y$.

If $y\notin S$, let $c=f(y)$ and $\tilde S_c=F_v\cap \partial \Omega_c$.
Consider the `inverse' of $\pi|_{\tilde S_c}: \tilde S_c\to S$,
$$\varphi_c\triangleq \iota\circ
(\pi|_{\tilde S_c})^{-1}:S\longrightarrow  X, \eqno{(4.1)}$$
where $\op{Im}(\varphi_c)=\tilde S_c$ and $\iota: \tilde S_c\hookrightarrow X$
denotes the inclusion; the reason for adding $\iota$ is because apriori $\tilde S_c$ may not be an Alexandrov space.

Observe that if $\varphi_c$ is a radial cone-neighborhood isometry, then by (1.10.1)
$\varphi_c$ induces a radial cone-neighborhood isometric tangent map,
$$\text{\rm D}\varphi_c:\Sigma_pS\to \Sigma_{y}X,\eqno{(4.2)}$$
where $\pi(y)=p$ ($y\in \varphi_c(p)$). Since $y$ is an interior point of
a ray $\subset F_v$, $\Sigma_yX\in \text{Alex}(1)$ is isometric to a spherical suspension over
the cross section $\Sigma_y\partial \Omega_c\in \text{Alex}(1)$.
Because $\tilde S_c=\varphi_c(S)\subset\partial\Omega_c$,
$\text{\rm D}\varphi_c(\Sigma_pS)\subset \Sigma_{y}\partial\Omega_c$.
(Note that although $\Sigma_{y}F_v$ is a spherical suspension over
$\Sigma_{y}\tilde S_c$, it is not clear if $\Sigma_{y}F_v\in \text{Alex}(1)$.)

If $y\in S$, by (3.1.2) $\Sigma_yF_v\cap\Uparrow_y^{\partial\Omega_c}
=\{\xi_1,\cdots,\xi_h\}$, $h\le k$. We will describe $\Sigma_yF_v$ as
follows: for fixed $\xi_i$, let $[\xi_i\Sigma_yS]'_{F_v}$
denote the union of the minimal geodesics of length $\frac \pi2$,
each of which is determined by a flat strip containing $\gamma_{y,\xi_i}$ and a
minimal geodesic in $S$ at $y$, and let $[\xi_i\Sigma_yS]_{F_v}$ denote the
closure of $[\xi_i\Sigma_yS]'_{F_v}$ in $\Sigma_yX$. Then
$$\Sigma_yF_v=\bigcup_{i=1}^h [\xi_i\Sigma_yS]_{F_v}.$$
Note that $[\xi_i\Sigma_yS]_{F_v}$ determines
a natural multi-valued map
$$\sigma_{\xi_i}:\Sigma_yS\to \Sigma_{\xi_i}(\Sigma_yX) \text{ by }
\eta\mapsto\Uparrow_{\xi_i,F_v}^{\eta},\eqno{(4.3)}$$
where $\Uparrow_{\xi_i,F_v}^{\eta}$ is the union of directions from $\xi_i$ to $\eta$ in $[\xi_i\Sigma_yS]_{F_v}$.

The main result in this section is:

\proclaim{Theorem 4.1}
\noindent {\rm (4.1.1)} The $\varphi_c$ in (4.1) is a radial cone-neighborhood isometry.

\noindent {\rm (4.1.2)} The tangent map in (4.2), $\text{\rm D}\varphi_c:\Sigma_pS\to \Sigma_{y}\partial\Omega_c\subset\Sigma_yX$, is a radial cone-neighborhood isometry
such that $\op{Im}(\text{\rm D}\varphi_c)=\Sigma_{y}\tilde S_c$. Moreover, if $|\xi\Sigma_{y}\tilde S_c|\geq\frac\pi2$ for some $\xi\in
\Sigma_{y}X$, then $|\xi\eta|=\frac\pi2\ \text{ for all } \eta\in\Sigma_{y}\tilde S_c$.

\noindent {\rm (4.1.3)} The multi-valued map $\sigma_{\xi_i}$ in (4.3) is a radial cone-neighborhood isometry.
\endproclaim

\subhead 4.1. Proof of (4.1.1) and its corollaries
\endsubhead

\demo{Proof of (4.1.1)}

Observe that $\varphi_c$ satisfies (1.9.1), because $|\varphi_c(q)|<|\varphi_c(p)|=k$
for all $p\in S_k(F_v)$ and $q\in S\setminus  S_k(F_v)$ (by (3.1.2) and Remark 3.4)
and $ S_k(F_v)$ is dense in $S$. Then we only need to show that $\varphi_c$ satisfies (1.9.2).

Note that for $p_1,p_2\in S$ and $\tilde p_i\in \varphi_c(p_i)$, $|\tilde p_1\tilde p_2|\geq |p_1p_2|$
because $\pi$ is distance non-increasing. By definition of $F_v$, given a
$[p_1p_2]$, there is a ray $\gamma_{w,p_2}$ such that $\{\gamma_{\uparrow_{p_1}^{\tilde p_1},p_1},[p_1p_2],\gamma_{w,p_2}\}$ bounds a flat strip $\Cal S\subset F_v$,
where $\gamma_{\uparrow_{p_1}^{\tilde p_1},p_1}$ denotes the ray from $p_1$ to
$\tilde p_1$. Note that $\Cal S\cap\tilde S_c$ is a minimal geodesic
parallel to $[p_1p_2]$ in $\Cal S$, i.e. $\varphi_c:[p_1p_2]\to \Cal S\cap\tilde S_c$ is an isometry.
It follows that $\varphi_c$ satisfies (1.9.2).
\hfill$\qed$
\enddemo

Applying (1.9.3) and (1.9.4) to the radial cone-neighborhood isometry $\varphi_c$ ((4.1.1)),
we conclude the following properties.

\proclaim{Corollary 4.2}

\noindent{\rm(4.2.1)} Let $p\in  S_k(F_v)$,
$\gamma_{w,p}\subset F_v$, and let $B_\epsilon(p)\subset S$ satisfy (1.9.3). Then
there is an isometric embedding
$B_\epsilon(p)\times [0,+\infty)\hookrightarrow F_v$
whose image contains $\gamma_{w,p}$.

\noindent{\rm(4.2.2)} Let $[pq]\subset S$ with $p\in S\setminus  S_k(F_v)$, $q\in S_k(F_v)$
close to $p$, and let $V_{[pq],\delta}$ satisfy (1.9.4). Suppose that $\gamma_{w,p}\subset F_v$
is parallel to $\gamma_{u,q}\subset F_v$ along $[pq]$. Then for small $\delta$, there
exists an isometric embedding $V_{[pq],\delta}\times
[0,+\infty)\hookrightarrow  F_v$ whose image contains $\gamma_{w,p}$ and $\gamma_{u,p}$.
\endproclaim

\subhead 4.2. Proof of (4.1.2)
\endsubhead

The former part of (4.1.2) follows from (4.1.1), (1.10.1) and the comments after (4.2).
And the latter part of (4.1.2) is an immediate corollary of Lemma 4.3 below (note that
we cannot apply (1.10.2) here because $\dim(\Sigma_{y}\partial\Omega_c)=\dim(S)+1$).

\proclaim{Lemma 4.3} Let $Y, \tilde Y\in \text{\rm Alex}(1)$ with
$\dim(\tilde Y)>\dim(Y)\geq1$ and $\partial Y=\emptyset$, and let $f:Y\rightarrow
\tilde Y$ be a radial cone-neighborhood isometry. If $|\xi f(Y)|\geq\frac\pi2$ for some $\xi\in \tilde Y$,  then
$|\xi\eta|=\frac\pi2\ \text{ for all } \eta\in f(Y)$.
\endproclaim

Note that if $f(Y)$ is locally convex in $\tilde Y$, then $\partial f(Y)=\emptyset$, and thus Lemma 4.3
follows from Lemma 1.3.

\demo{Proof} We proceed by induction on $\dim(Y)$.
When $\dim(Y)=1$, $Y$ is a circle. By (1.9.2) (and (1.9.3)), it is easy to see that $f(Y)$ is a union of several
disjoint locally convex circles.  Then the desired property follows by
applying Lemma 1.3 to $f(Y)$.

Assume that $\dim(Y)>1$. Let $\xi\in\tilde Y$ with $|\xi f(Y)|\ge \frac \pi2$.
Claim 1: there is $\eta_0\in f(Y)$ such that
$|\xi\eta_0|=|\xi f(Y)|$. By a standard open and closed argument, it reduces to show
that for $\eta\in f(Y)$ close to $\eta_0$, $|\xi \eta|=\frac \pi2$. Claim 2:
there is $[\xi\eta_0]$ and $[\eta_0\eta]$ such that $|\uparrow_{\eta_0}^\xi\uparrow_{\eta_0}^{\eta}|=\frac\pi2$.
By (1.7.2), ``$|\xi \eta|\ge |\xi \eta_0|\ge \frac \pi2$'' and ``$|\uparrow_{\eta_0}^\xi\uparrow_{\eta_0}^{\eta}|=\frac\pi2$''
implies that $$|\xi\eta_0|=|\xi\eta|=\frac\pi2.$$
Hence, in the rest we only need to verify Claim 1 and 2.

Note that Claim 1 follows if we show that $f(Y)$ is closed and thus compact.
Let $\tilde y$ be a limit point of $f(Y)$, i.e. there are $\tilde y_i\in f(Y)$ such that
$\tilde y_i\to \tilde y$. Let $y_i\in Y$ such that $f(y_i)\ni\tilde y_i$. Passing to a subsequence,
we can assume that $y_i\to y$. By (1.9.2), it is easy to see that $\tilde y\in f(y)$.

In order to verify Claim 2, we
Let $\bar\eta_0\in Y$ such that $\eta_0\in f(\bar\eta_0)$.
Since $f$ is a radial cone-neighborhood isometry, for $\eta\in f(Y)$ close to $\eta_0$
there is $[\eta_0\eta]\subset f(Y)$
with $\uparrow_{\eta_0}^{\eta}\in\text{D}f(\Sigma_{\bar\eta_0}Y)$ (see (1.9.3-4)),
where $\text{D}f:\Sigma_{\bar\eta_0}Y\to\Sigma_{\eta_0}\tilde Y$ is a tangent map
of $f$ (see (1.10.1)). On the other hand, because $|\xi\eta_0|=|\xi f(Y)|$,
for any $[\eta_0\xi]$ we have that $|\uparrow_{\eta_0}^\xi\text{D}f(\Sigma_{\bar\eta_0}Y)|\geq\frac\pi2$.
By (1.10.1), $\text{D}f$ is also a radial cone-neighborhood isometry, so by the
inductive assumption
$$|\uparrow_{\eta_0}^\xi\zeta|=\frac\pi2 \text{ for all } \zeta\in\text{D}f(\Sigma_{\bar\eta_0}Y).$$
Hence, $|\uparrow_{\eta_0}^\xi\uparrow_{\eta_0}^{\eta}|=\frac\pi2$, i.e. Claim 2 holds.
\hfill$\qed$
\enddemo

\subhead 4.3. Proof of (4.1.3)
\endsubhead

Note that if $y\in  S_k(F_v)$, then by (4.2.1)
$[\xi_i\Sigma_yS]_{F_v}=\{\xi_i\}*\Sigma_yS\ (\subset\Sigma_yX)$,
so each $\sigma_{\xi_i}$
is an isometric embedding, and thus a radial cone-neighborhood isometry.

If $y\in S\setminus  S_k(F_v)$, we need to show that $\sigma_{\xi_i}$ satisfies (1.9.1) and (1.9.2).

Note that for any $B_\epsilon(y)\subset S$, $\bigcup\limits_{y'\in B_\epsilon(y)\cap S_k(F_v)}\Uparrow_{y}^{y'}$
is dense in $\Sigma_yS$ because $ S_k(F_v)$ is dense in $S$. By (3.11) and (3.12), for small $\epsilon$,
there is $l\leq k$ such that for all $\eta\in\Sigma_yS$ and $\uparrow_{y}^{y'}\in\Uparrow_{y}^{y'}$
$$|\Uparrow_{\xi_i,F_v}^{\eta}|\leq |\Uparrow_{\xi_i,F_v}^{\uparrow_{y}^{y'}}|=l.$$
That is, $\sigma_{\xi_i}$ satisfies (1.9.1).

On the other hand, by (1.7.2) it is clear that $|\sigma_{\xi_i}(\eta)\sigma_{\xi_i}(\eta')|\geq|\eta\eta'|$ for all $\eta,\eta'\in\Sigma_yS$.
Then Lemma 4.4 below implies that $\sigma_{\xi_i}$ satisfies (1.9.2).
\hfill$\qed$

\proclaim{Lemma 4.4} For any $[\eta\eta']\subset\Sigma_yS$ and $[\xi_i\eta]\subset[\xi_i\Sigma_yS]_{F_v}$, there is an isometric embedding $\{\xi_i\}*[\eta\eta']\hookrightarrow [\xi_i\Sigma_yS]_{F_v}$ such that $\{\xi_i\}*\{\eta\}=[\xi_i\eta]$.
\endproclaim

\demo{Proof} Since $\Sigma_y'S$ is dense in $\Sigma_yS$, we
can assume that $\eta,\eta'\in\Sigma_y'S$, i.e. there is $[yq]$ and $[yq']$ (in $S$)
such that $\eta=\uparrow_y^q$ and $\eta'=\uparrow_y^{q'}$. We select two sequences
$q_j\in[yq], q_j'\in[yq']$ with
$|q_jy|=|q_j'y|\to 0$  as $j\to\infty$.
Note that there is $[q_jq_j']\ (\subset S)$ such that $\bigcup\limits_{p\in[q_jq_j']}\Uparrow_{y}^{p}$
converges to $[\eta\eta']$ (cf. [BGP]).

Let $\gamma_{w_j, q_j}\subset F_v$ be the ray
such that $[\xi_i\eta]$ is determined by the flat strip $\Cal S_j$
bounded by $\{\gamma_{\xi_i,y},[yq_j], \gamma_{w_j, q_j}$\}. By the construction
of $F_v$, there is $\gamma_{w_j', q_j'}\subset F_v$ such that
$\{\gamma_{w_j,q_j},[q_jq_j'], \gamma_{w_j', q_j'}\}$
bounds a flat strip in $F_v$.
And by (3.10), $\{\gamma_{w_j',q_j'}, [yq_j'], \gamma_{\xi_i, y}\}$ also
bounds a flat strip $\Cal{S}_j'\subset F_v$.  Note that
all $\Cal S_j$ lie in a flat strip in $F_v$ over $[yq]$.
Hence, if all $\Cal S_j'$ also lie in a flat strip in $F_v$ over $[yq']$,
then $[\xi_i\eta], [\eta\eta']$ and $[\xi_i\eta']$ determined by $\Cal S_j'$
bound a convex spherical surface in $[\xi_i\Sigma_yS]_{F_v}$, i.e. Lemma 4.4. holds.
In fact, by (3.1.2) there are only a finite number of flat strips over $[yq']$ in $F_v$, so
passing to a subsequence we can assume that $\Cal S_j'$ lies in a same flat strip over $[yq']$ in $F_v$.
\hfill$\qed$
\enddemo

%%%%%%%%%%%%%%%%%%%%%%%%%%%%%%%% section 6 %%%%%%%%%%%%%%%%%%%%%

\vskip4mm

\head 5. Proof of (3.1.3)
\endhead

\vskip4mm

Our proof of (3.1.3) is based on the following lemma.

\proclaim{Lemma 5.1}  Let $Y\in\text{\rm Alex}^{n}(1)$ with
$n\geq2$, and let $Z$ be a locally convex subset in $Y$ with  $\dim(Z)=n-2$
and $\partial Z=\emptyset$. Suppose that $|pZ|\geq\frac\pi2$, $p\in Y$.
If $f_i:Z\to \Sigma_pY$ by $z\mapsto (\Uparrow_i)_p^z$ ($i=1,2$), where
$(\Uparrow_i)_p^z\subset \Uparrow_p^z$, are two
radial cone-neighborhood isometries, then $f_1=f_2$.
\endproclaim

\demo {Proof of (3.1.3) by assuming Lemma 5.1}

We argue by contradiction. If $F_v\ne F_w$ and $F_v\cap F_w\supsetneq S$,
$F_v\cap F_w$ contains a ray $\gamma_{\xi,q}$ with $q\in S$ and
$\xi\in\Uparrow_{q}^{\partial\Omega_c}$.
Note that $F_v$ and $F_w$ determine  two radial
cone-neighborhood isometries from $\Sigma_{q}S$ to
$\Sigma_{\xi}(\Sigma_{q}X)$ (see (4.1.3)).
Moreover, note that
$$\xi\in(\Sigma_{q}S)^\perp \text{ and }
\dim(\Sigma_{q}X)=\dim(\Sigma_{q}S)+2.$$ It
therefore follows from Lemma 5.1 that such two radial
cone-neighborhood isometries coincide, which implies that
$F_v=F_w$, a contradiction. \hfill$\qed$
\enddemo

Based on (3.1.3), one may observe a `radial' total convexity of $F_v$.

\proclaim{Corollary 5.2} Let $\alpha_i(t)|_{ [0,+\infty)}$
($i=1,2$) be two rays in $F_v$ with $\alpha_i(0)\in S$. If
$\alpha_1(t)|_{[0,+\infty)}$, $\alpha_2(t)|_{[0,+\infty)}$
and some $[\alpha_1(0)\alpha_2(0)]$ bounds a flat strip, then  for all $t_i\in [0,+\infty)$  any
$[\alpha_1(t_1)\alpha_2(t_2)]$ belongs
to $F_v$.
\endproclaim

\demo {Proof} For convenience, we let $x_i$ denote $\alpha_i(t_i)$.
Without loss of generality, we assume that
$|x_1\alpha_1(0)|\geq|x_2\alpha_2(0)|$ and $|x_1\alpha_1(0)|>0$.
Since $\alpha_1(t)|_{[0,+\infty)}$,
$\alpha_2(t)|_{[0,+\infty)}$ and some $[\alpha_1(0)\alpha_2(0)]$
bounds a flat strip, by (1.7.2) it is not hard to see that
$$|\alpha_1^-(t_1)\uparrow_{x_1}^{x_2}|=
\arccos\frac{f(x_2)-f(x_1)}{|x_1x_2|},\eqno{(5.1)}$$ where $f$ is
the Busemann function (note that
$f(x_2)-f(x_1)=|x_1\alpha_1(0)|-|x_2\alpha_2(0)|$ by Lemma 1.2).
On the other hand, since $f$ is concave (see Lemma 1.1) and
$\nabla_{\alpha_i(t)}f=\alpha_i^-(t)$ (by Lemma 1.2), we have
that ([Pet1])
$$\text{d}f_{x_1}(\uparrow_{x_1}^{x_2})\leq
\cos|\alpha_1^-(t_1)\uparrow_{x_1}^{x_2}|.$$ Then by
the concavity of $f$ and (5.1), $f|_{[x_1x_2]}$ must be a linear
function, i.e.
$$f(y)=f(x_1)+|x_1y|\cos|\alpha_1^-(t_1)\uparrow_{x_1}^{x_2}|
\text{ for any } y\in [x_1x_2].\eqno{(5.2)}$$

Furthermore, by (1.7.3), (5.1) implies that there is an
$[x_2\alpha_1(0)]$ such that $[x_1x_2]$, $[x_1\alpha_1(0)]$ and
$[x_2\alpha_1(0)]$ bounds a convex flat surface $\Cal E$. For any
$y\in [x_1x_2]$, we let $y'\in [x_2\alpha_1(0)]$ be the point such
that $[yy']$ in $\Cal E$ is parallel to $[x_1\alpha_1(0)]$. Similarly,
we can calculate $f(y')$ as $f(y)$ in (5.2), and it is easy to see that
$$f(y')=f(y)+|yy'|.$$
Then by Lemma 1.2, we can conclude that $[y'y]$ belongs to a ray
starting from $y'$ with $\uparrow_{y'}^{y}\in
\Uparrow_{y'}^{\partial\Omega_{c}}$ for $c<f(y')$.

Similarly, we can find a minimal geodesic $[\alpha_1(0)\alpha_2(0)]'$
(which may be $[\alpha_1(0)\alpha_2(0)]$) such that
$[x_2\alpha_1(0)]$, $[x_2\alpha_2(0)]$ and
$[\alpha_1(0)\alpha_2(0)]'$ bounds a convex flat surface. And there
is a ray, which starts from some $y''\in [\alpha_1(0)\alpha_2(0)]'$
and is perpendicular to $\partial\Omega_{c}$, such that both $y'$
and $y$ lie in it. In other words, $[x_1x_2]$ belongs to an
$F_w\subset\Cal F$. By (3.1.3), it has to hold  that
$F_w= F_v$, i.e., $[x_1x_2]\subset  F_v$.
\hfill$\qed$
\enddemo

As a preparation for Lemma 5.1, we give the following partial analogy of
Frankel's Theorem in Riemannian geometry (cf. [Pet2]).

\proclaim{Lemma 5.3}  Let $Y\in\text{\rm Alex}^{n}(1)$ with $n\geq2$
and $Y_i\in\text{\rm Alex}^{n_i}(1)$ with $n_1=n-1$, $n_2\geq1$ and
$\partial Y_i=\emptyset$. If there are two radial cone-neighborhood
isometries $f_i:Y_i\to Y$ (in  particular, each $f_i$ is an isometric
embedding), then $f_1(Y_1)\cap f_2(Y_2)\neq\emptyset$.
\endproclaim

\demo {Proof of Lemma 5.1 by assuming Lemma 5.3}

We give the proof by induction on $n$ starting with
$n=2$.

If $n=2$, then $\dim (Z)=0$. In this case, that $Z$ is locally convex and
$\partial Z=\emptyset$ means that $Z$ consists of two points with
distance $\pi$, so the conclusion is clear.

Assume that $n>2$. By Lemma 5.3,
$f_1(Z)\cap f_2(Z)\neq\emptyset$, i.e., there is $[pz]$ with $z\in
Z$ such that $\uparrow_p^z\in f_1(Z)\cap f_2(Z)$. By (1.10.1), the
tangent maps $\text{D}f_i:\Sigma_zZ\to\Sigma_{\uparrow_p^z}Y$ are
also radial cone-neighborhood isometries. By the inductive
assumption, we have that $\text{D}f_1=\text{D}f_2$, which implies
that $f_1=f_2$. \hfill$\qed$
\enddemo

In the rest of this section, we will complete the proof of (3.1.3) by
verifying Lemma 5.3.

\demo {Proof of Lemma 5.3}

For convenience, we let $A\triangleq f_1(Y_1)$ and
$B\triangleq f_2(Y_2)$. If $A\cap B=\emptyset$, then there is $a\in
A$ and $b\in B$ such that $|ab|=|AB|>0$. We claim that
$$|ab|<\frac\pi2.\eqno{(5.3)}$$
If the claim fails, then by Lemma 4.3 $B\subseteq
A^{=\frac\pi2}=A^{\geq\frac\pi2}$ (and $A\subseteq B^{=\frac\pi2}$).
By (1.8.1) and (1.9.3), we have that
$\dim(Y_1)+\dim(Y_2)\leq n-1$, which contradicts
``$\dim(Y_1)+\dim(Y_2)\geq n$''.

With (5.3), we will show the existence, $a'\in A$ and $b'\in B$, such that
$|a'b'|<|ab|$, a contradiction.

By (1.10.1), $\text{\rm D} f_1:\Sigma_{\hat a} Y_1 \to \Sigma_aY$ is
also a radial cone-neighborhood isometry, where $\hat a\in Y_1$ with
$f_1(\hat a)\ni a$. Again by Lemma 4.3, for any fixed $[ab]$, we have
that $|\uparrow_a^b\eta|=\frac\pi2$ for all $\eta\in\Sigma_aA$. And
note that $\dim(A)=n_1=n-1$ (see (1.9.3)). It then is not hard to see that $g:\Sigma_aA\to
\Sigma_{\uparrow_a^b}(\Sigma_aY)$ by
$\eta\mapsto\Uparrow_{\uparrow_a^b}^\eta$ is an at most 2-valued
map, which implies that the following composition of maps,
$$g\circ\text{\rm D} f_1:\Sigma_{\hat a} Y_1 \to \Sigma_aA \to
\Sigma_{\uparrow_a^b}(\Sigma_aY) \text{ is a radial
cone-neighborhood isometry}.\eqno{(5.4)}$$ Since $\dim(Y_1)=n-1$ and $\partial
Y_1=\emptyset$, $g\circ\text{\rm D} f_1$ has to be surjective (by
(1.10.2)). It follows that for any $\xi\in\Sigma_aY$ with
$|\xi\uparrow_a^b|\leq\frac\pi2$ there exists $\eta\in\Sigma_aA$
such that
$$|\uparrow_a^b\xi|+|\xi\eta|=\frac\pi2\ (=|\uparrow_a^b\eta|). \eqno{(5.5)}$$
Moreover, for any other $\eta'\in \Sigma_aA$, we have that
$$\cos|\xi\eta'|=\cos|\xi\eta|\cos|\eta\eta'|. \eqno{(5.6)}$$

On the other hand, since $f_2$ is a radial cone-neighborhood
isometry, there is $[bb_0]$ (in $Y$) which lies in $B$.
Similarly, we can conclude that
$$|\uparrow_b^a\uparrow_b^{b_0}|=\frac\pi2.\eqno{(5.7)}$$
Let $b_i\in [bb_0]$ such that $b_i\to b$ as $i\to \infty$. Without loss
of generality, we may assume that $[ab_i]\to [ab]$. By (5.5), there is $\eta_i\in
\Sigma_aA$ such that
$$|\uparrow_a^b\uparrow_a^{b_i}|+|\uparrow_a^{b_i}\eta_i|=\frac\pi2.\eqno{(5.8)}$$
Passing to a subsequence, $\eta_i\to\eta_0$ as
$i\to\infty$. Note that we can find $[aa_0]\subset A$ such that
$$|\uparrow_a^{a_0}\eta_0|<<|ab|.\eqno{(5.9)}$$
Put $|bb_i|=\epsilon_i$, and for $i$ large let $a_i\in [aa_0]$ such
that $|aa_i|=\epsilon_i$. Let
$\triangle\tilde a\tilde b\tilde b_i$ be a comparison triangle (in
$\Bbb S^2$) of $\triangle abb_i$. Then we derive the following:
$$\align &\cos|a_ib_i|\\
\overset{\text{by (1.7.2)}}\to\geq&\cos|ab_i|\cos\epsilon_i+
\sin|ab_i|\sin\epsilon_i\cos|\uparrow_a^{b_i}\uparrow_a^{a_0}|\\
\overset{\text{by Lemma 4.3, (1.7.2)}}\to\geq&\cos|ab|\cos^2\epsilon_i+
\sin|ab_i|\sin\epsilon_i\cos|\uparrow_a^{b_i}\uparrow_a^{a_0}|\\
\overset{\text{by (5.6)}}\to=&\cos|ab|\cos^2\epsilon_i+
\sin|ab_i|\sin\epsilon_i\cos|\uparrow_a^{b_i}\eta_i|\cos|\eta_i\uparrow_a^{a_0}|\\
\overset{\text{by (5.8)}}\to=&\cos|ab|(\cos^2\epsilon_i+
\sin|ab_i|\sin\epsilon_i\sin|\uparrow_a^{b_i}\uparrow_a^{b}|
\frac{\cos|\eta_i\uparrow_a^{a_0}|}{\cos|ab|})\\
\overset{\text{by (1.7.2), (5.3),
(5.9)}}\to\geq&\cos|ab|(\cos^2\epsilon_i+
\sin|ab_i|\sin\epsilon_i\sin|\uparrow_{\tilde a}^{\tilde
b_i}\uparrow_{\tilde a}^{\tilde b}|
\frac{\cos|\eta_i\uparrow_a^{a_0}|}{\cos|ab|})\\
\overset{\text{by law of sine}}\to=&\cos|ab|(\cos^2\epsilon_i+
\sin^2\epsilon_i\sin|\uparrow_{\tilde b}^{\tilde
b_i}\uparrow_{\tilde b}^{\tilde a}|
\frac{\cos|\eta_i\uparrow_a^{a_0}|}{\cos|ab|})\\
\overset{\text{by (5.7), (5.9) and } i\to\infty}\to > &\cos|ab|.
\endalign$$
We now specify $a'=a_i$ and $b'=b_i$ for large $i$.
\hfill$\qed$
\enddemo

%%%%%%%%%%%%%%%%%%%%%%%%%%%%%%%% section 7 %%%%%%%%%%%%%%%%%%%%%

\vskip4mm

\head 6. Proof of (0.6.1)
\endhead

\vskip4mm

Recall that Lemma 2.1 provides a criterion for $d(\cdot ,Z)$ to be concave. Let
$x$ be an arbitrary point in $X\setminus Z$,
$\alpha(t)|_{[0, \epsilon)}$, $y$, $[xy]$ and $[y\alpha(t)]$ as in Lemma 2.1.
By specifying $Z=\Cal F$, the concavity of $d(\cdot,\Cal F)$ follows from
(2.1.1) and (2.1.3). Our verification of (2.1.1) and (2.1.3) is technical
and tedious, and is divided into two cases: $y\notin S$ or $y\in S$.

\demo{Proof of (0.6.1) for $y\not\in S$, i.e., $y$ lies in some $F_v\setminus S$}

We first claim that $(\uparrow_y^x)^\perp=\Sigma_y F_v$
and $\Uparrow_{\uparrow_y^x}^{(\uparrow_y^x)^\perp}=\Sigma_{\uparrow_y^x}(\Sigma_y X)$
(i.e., $B_{\frac\pi2}(\uparrow_y^x)=[\uparrow_y^x(\uparrow_y^x)^\perp]$);
and for $\eta,\eta'\in(\uparrow_y^x)^\perp$ with $\eta'$ close to $\eta$, any $[\uparrow_y^x\eta]$ is perpendicular to an $[\eta\eta']$ (e.g.,
$[\uparrow_y^x(\uparrow_y^x)^\perp]=\{\uparrow_y^x\}*(\uparrow_y^x)^\perp$).
Note that (2.1.1) and (2.1.3.1) follows from the former part of the claim,
and (2.1.3.2) follows from the latter part and (1.7.3).

We now verify the above claim.
Let $\beta(t)|_{[0,+\infty)}\subset F_v$
with $\beta(0)\in S$ be the ray such that $y=\beta(t_0)$ with $t_0>0$. Let $c\triangleq f(y)$ and $\tilde S_c\triangleq
F_v\cap\partial\Omega_c$. Note that
$\Sigma_yX=\{\beta^+(t_0),\beta^{-}(t_0)\}*\Sigma_y\partial\Omega_c$
and $\Sigma_yF_v=\{\beta^+(t_0),\beta^{-}(t_0)\}*\Sigma_y\tilde S_c$. Hence, it suffices to show that $\uparrow_y^x\in \Sigma_y\partial\Omega_c$, $|\uparrow_y^x\eta|=\frac\pi2$ for any $\eta\in \Sigma_y \tilde S_c$, $\Uparrow_{\uparrow_y^x}^{\Sigma_y \tilde S_c}=\Sigma_{\uparrow_y^x}(\Sigma_y \partial\Omega_c)$,
and any $[\uparrow_y^x\eta]$ is perpendicular to any $[\eta\eta']\subset\Sigma_y \tilde S_c$ with $\eta'$ close to $\eta$.

Since $|xy|=|x\Cal F|$, by the first variation formula, we derive that
$|\uparrow_y^x\Sigma_yF_v|\geq\frac\pi2$. Hence, $\uparrow_y^x\in \Sigma_y\partial\Omega_c$, and by (4.1.2)
$|\uparrow_y^x\eta|=\frac\pi2$ for any $\eta\in \Sigma_y\tilde S_c$.
Consider the multi-valued map (see (4.1.2)),
$$\psi\circ\text{D}\varphi_c:\Sigma_{\beta(0)}S\overset{\text{D}\varphi_c}\to\longrightarrow
\Sigma_y\tilde S_c\overset{\psi}\to\longrightarrow
\Sigma_{\uparrow_y^x}(\Sigma_y\partial\Omega_c), \eqno{(6.1)}$$
where
$\psi(\eta)=\Uparrow_{\uparrow_y^x}^{\eta}$ for any $\eta\in
\Sigma_y\tilde S_c$. Note that ``$\Uparrow_{\uparrow_y^x}^{\Sigma_y \tilde S_c}=\Sigma_{\uparrow_y^x}(\Sigma_y \partial\Omega_c)$'' is equivalent to ``$\psi$ is surjective''.
By (4.1.2) and (1.7.3), similar to (5.4) we conclude that $\psi\circ\text{D}\varphi_c$ is
a radial cone-neighborhood isometry. Because $\dim(\Sigma_{\beta(0)}S)=\dim(\Sigma_{\uparrow_y^x}(\Sigma_y\partial\Omega_c))$
and $\partial\Sigma_{\beta(0)}S=\emptyset$, by (1.10.2)
$\psi\circ\text{D}\varphi_c$  is surjective, so is $\psi$ surjective.

Furthermore, by (1.10.1) $\Sigma_\eta(\Sigma_y \tilde S_c)$ is the image
of a tangent map of $\text{D}\varphi_c$; and
$|\uparrow_{\eta}^{\uparrow_y^x}\Sigma_\eta(\Sigma_y \tilde S_c)|\geq\frac\pi2$ because
$|\uparrow_y^x\eta|=|\uparrow_y^x\Sigma_y \tilde S_c|=\frac\pi2$.
By Lemma 4.3, $[\uparrow_y^x\eta]$ is perpendicular to $[\eta\eta']$
with $\eta'$ close to $\eta$.

Finally, by the definition of $F_v$, for any $u\in\Sigma_y F_v$ and small $t$,
there is $[yy']\subset F_v$ such that $|yy'|=t$ and $\uparrow_y^{y'}\to u$ as $t\to 0$ (cf. Corollary 5.2).
Note that this guarantees (2.1.3.3).
\hfill$\qed$
\enddemo

For $y\in S$, similar to the claim in the above proof, we need the following lemma.

\proclaim{Lemma 6.1} If $y\in S$, then the following holds.

\noindent {\rm (6.1.1)} $\{\uparrow_{y}^x,\Uparrow_{y}^{\partial\Omega_c}\}\subset(\Sigma_yS)^\perp$ and
$|\uparrow_{y}^x\Uparrow_{y}^{\partial\Omega_c}|=\frac\pi2$, and thus $\dim((\Sigma_yS)^\perp)=1$.
As a result, $(\uparrow_y^x)^\perp\cap \Uparrow_{y}^{\partial\Omega_c}=\{\xi_1,\xi_2\}$,
where $\xi_1$ may be equal to $\xi_2$.

\noindent {\rm (6.1.2)} $\Uparrow_{\uparrow_y^x}^{(\uparrow_y^x)^\perp}=\Sigma_{\uparrow_y^x}(\Sigma_y X)$
(i.e., $B_{\frac\pi2}(\uparrow_y^x)=[\uparrow_y^x(\uparrow_y^x)^\perp]$);
and for $\zeta,\zeta'\in(\uparrow_y^x)^\perp$ with $\zeta'$ close to $\zeta$, any $[\uparrow_y^x\zeta]$ is perpendicular to an $[\zeta\zeta']$.

\noindent {\rm (6.1.3)} $(\uparrow_y^x)^\perp\subset \Sigma_y\Cal F$; in detail,  $(\uparrow_y^x)^\perp=\bigcup_{i=1}^2[\xi_i\Sigma_yS]_{F_{v_i}}$, where $F_{v_i}$ contains $\gamma_{\xi_i,y}$
(see (4.3) for $[\xi_i\Sigma_yS]_{F_{v_i}}$).
\endproclaim

\demo{Proof of (0.6.1) for $y\in S$ by assuming Lemma 6.1}

It is easy to see that (2.1.1) and (2.1.3.1) follows from the first
property in (6.1.2) and (6.1.3), and (2.1.3.2) follows from the second property in
(6.1.2) and (1.7.3). Moreover, since $(\uparrow_y^x)^\perp=\bigcup_{i=1}^2[\xi_i\Sigma_yS]_{F_{v_i}}$ ((6.1.3)), for any $u\in (\uparrow_y^x)^\perp$ and small $t$, there is a flat strip
$\Cal S_j\subset F_{v_i}$ with $\gamma_{\xi_i,y}\subset \Cal S_j$ ($i=1$ or 2)
such that we can select $[yy_j]\subset \Cal S_j$ with $\uparrow_y^{y_j}\in(\uparrow_y^x)^\perp$ so that
$|yy_j|=t\to0$ and $\uparrow_y^{y_j}\to u$ as $j\to\infty$.
This suffices to see (2.1.3.3).
\hfill$\qed$
\enddemo

The main effort in the rest of Section 6 is to verify Lemma 6.1.

\demo{Proof of (6.1.1)}

Note that $\Uparrow_{y}^{\partial\Omega_c}\subset(\Sigma_yS)^\perp$ (see (1.1)).
On the other hand, since $|xy|=|x\Cal F|$, $|\uparrow_y^x\Sigma_yS|\geq\frac\pi2$ and
$|\uparrow_y^x\Uparrow_{y}^{\partial\Omega_c}|\geq\frac\pi2$.
By Lemma 1.3, ``$|\uparrow_y^x\Sigma_yS|\geq\frac\pi2$'' implies that $\uparrow_y^x\in(\Sigma_yS)^\perp$.
And because the concave Busemann function
$f$ achieves the maximum  on $S$, we have
$$|\Uparrow_y^{\partial\Omega_c}v|\leq\frac\pi2\ (\text{i.e. } df|_y(v)\leq0).\eqno{(6.2)}$$
Hence, it follows that $$|\uparrow_y^x\Uparrow_{y}^{\partial\Omega_c}|=\frac\pi2.\eqno{(6.3)}$$
By (1.8.1), we notice that $\dim((\Sigma_yS)^\perp)\leq1$
(note that $(\Sigma_yS)^\perp$ is convex in $\Sigma_yX$). And by the convexity of $(\Sigma_yS)^\perp$,
$\dim((\Sigma_yS)^\perp)=1$ because $\{\uparrow_y^x,\Uparrow_{y}^{\partial\Omega_c}\}\subset(\Sigma_yS)^\perp$.
As a result, $(\uparrow_y^x)^\perp\cap \Uparrow_{y}^{\partial\Omega_c}$ contains at most two points.
\hfill$\qed$
\enddemo

In proving (6.1.2), the following technical result in [RW] is required.
In fact, Lemma 6.2 plays a role in the proof of Theorem 1.8.

\proclaim{Lemma 6.2 ([RW])} Let $Y\in \text{\rm Alex}^n(1)$, and
let $Y_1, Y_2$ be two compact convex subsets in $Y$
with $\dim(Y_1)+\dim(Y_2)=n-1$ and
$|y_1y_2|=\frac \pi2$ for any $y_i\in Y_i$. Then

\vskip1mm

\noindent {\rm (6.2.1)} Given $p_1\in
Y_1^\circ$, there is $m$ such that $\lambda_{p_1p_2}\leq m$ for
any $p_2\in Y_2^\circ$, where
$\lambda_{p_1p_2}$ denotes the number of minimal geodesics between
$p_1$ and $p_2$; and $\{p_2\in Y_2^\circ|\
\lambda_{p_1p_2}=m\}$ is open and dense in $Y_2^\circ$ (in
particular, it is just $Y_2^\circ$ if $n_2=1$). Moreover, for any $p_2\in Y_2^\circ$
with $\lambda_{p_1p_2}=m$ and any $[p_1p_2]$, there is a
neighborhood $U$ of $p_2$ in $Y_2^\circ$ such that $\{p_1\}* U$ can be embedded
isometrically into $Y$ around $[p_1p_2]$.

\vskip1mm

\noindent {\rm (6.2.2)} If  $\partial
Y_2=\emptyset$, then for any $p_1\in Y_1^\circ$ or for
$p_1\in\partial Y_1$ with
$(\Sigma_{p_1}Y_1)^{\geq\frac\pi2}=(\Sigma_{p_1}Y_1)^\perp$ we have that
$(\Sigma_{p_1}Y_1)^\perp=\Uparrow_{p_1}^{Y_2}$ and
$\partial(\Sigma_{p_1}Y_1)^\perp=\emptyset$. Moreover, the
multi-valued map $\sigma:Y_2\to (\Sigma_{p_1}Y_1)^\perp$ defined
by $p_2\mapsto\Uparrow_{p_1}^{p_2}$ is a radial cone-neighborhood
isometry; as a corollary,
$\sigma^{-1}$ is a metric cover$\ $
\footnote{In the paper, a metric cover means a locally isomeric cover.} if $n_2=1$.

\vskip1mm

\noindent {\rm (6.2.3)} If $Y_1=\{p_1\}$ and if
$\partial Y_2=\emptyset$, then either $\{p_1\}*Y_2$ can be embedded
isometrically into $Y$, or $Y\overset{\text{iso}}\to\cong
[\{p_1,\tilde p_1\}*\Sigma_{p_1}Y]/\Bbb Z_2$, where
$|p_1\tilde p_1|=\pi$ and the $\Bbb Z_2$-action maps $p_1$ to $\tilde p_1$.
\endproclaim

\demo{Proof of (6.1.2)}

By (6.1.1), $\dim((\Sigma_yS)^\perp)=1$, i.e. $(\Sigma_yS)^\perp$ is isometric to a circle
$S^1$ or an $[ab]$ with diameter $\leq\pi$.

If $(\Sigma_yS)^\perp=S^1$, then $\op{diam}(S^1)\ge \frac\pi2$ (note that $|\uparrow_y^x\xi_i|=\frac\pi2$ in $(\Sigma_yS)^\perp$), and thus
by (1.8.2), $\Sigma_yX=(\Sigma_yS)*S^1$, or
$\Sigma_yX=(\widetilde{\Sigma_yS}*\widetilde{S^1})/\Bbb Z_2$ with
$\widetilde{\Sigma_yS}/\Bbb Z_2=\Sigma_yS$ and $\widetilde{S^1}/\Bbb Z_2=S^1$
(only when $\xi_1=\xi_2$), which implies (6.1.2).
Consequently, in the rest of the proof we may assume that
$(\Sigma_yS)^\perp=[ab]$. We will divide the proof in two cases:
$\uparrow_y^x\in[ab]^\circ$ or otherwise.

\vskip1mm

\noindent Case 1: $\uparrow_y^x\in[ab]^\circ$.

\vskip1mm

By (6.2) and (6.3), it has to hold that
$(\uparrow_y^x)^\perp\cap \Uparrow_{y}^{\partial\Omega_c}=\Uparrow_y^{\partial\Omega_c}=\{a,b\}$
with $|ab|=\pi$. Hence, $\Sigma_yX=\{a,b\}*A$ for some
$A\in\text{Alex}^{n-2}(1)$. Note that
$\{\uparrow_y^x,\Sigma_yS\}\subset A$, so
$(\uparrow_y^x)^\perp=\{a,b\}*\Sigma_yS$ with
$\partial(\uparrow_y^x)^\perp=\emptyset$.
By (6.2.3), either $\{\uparrow_y^x\}*(\uparrow_y^x)^\perp$ can be
embedded isometrically into $\Sigma_yX$ (with full dimension) or
$\Sigma_yX$ is isometric to $(\{\uparrow_y^x,
\widetilde{\uparrow_y^x}\}*\widetilde{(\uparrow_y^x)^\perp})/\Bbb Z_2$
with $\widetilde{(\uparrow_y^x)^\perp}/\Bbb Z_2=(\uparrow_y^x)^\perp$, which implies (6.1.2).

\vskip1mm

\noindent Case 2: $\uparrow_y^x\notin[ab]^\circ$, i.e., $\uparrow_y^x=a$ or
$b$, say $a$.

\vskip1mm

We first claim that
$$[a\Sigma_yS]=\{a\}*\Sigma_yS \text{\ \ or\ \ } (\{a,\tilde
a\}*\widetilde{\Sigma_yS})/\Bbb Z_2 \text{ with }
\widetilde{\Sigma_yS}/\Bbb Z_2=\Sigma_yS.$$
Let $\xi\in [ab]$ with $|a\xi|=\frac\pi2$.
To see the claim, it suffices to show that
$$[a\Sigma_yS]\subseteq\{\xi\}^{=\frac\pi2} \text{ and }
[a\Sigma_yS] \text{ is totally convex in } \Sigma_yX.\eqno{(6.4)}$$
In fact, by (1.8.1), (6.4) implies that $\dim([a\Sigma_yS])=\dim(\Sigma_yX)-1=\dim(\Sigma_yS)+1$.
Then applying (6.2.3) to $[a\Sigma_yS]\in\text{Alex}(1)$ with $\partial(\Sigma_yS)=\emptyset$,
we can conclude the claim. We now verify (6.4).
Note that $(\uparrow_y^x)^\perp\cap \Uparrow_{y}^{\partial\Omega_c}=\{\xi\}$, so
$\Uparrow_y^{\partial\Omega_c}\subseteq[\xi b]$. Then by the first
variation formula, (6.2) implies that
$$|\uparrow_a^b\nu|\leq\frac\pi2 \text{ for all } \nu\in\Sigma_a(\Sigma_yX)\eqno{(6.5)}$$
(note that there is a unique minimal geodesic between $a$ and $b$).
It follows that
$(\uparrow_a^b)^\perp=(\uparrow_a^b)^{\geq\frac\pi2}$.
Then applying (6.2.2) to $\Sigma_yX$ (with $[ab]=Y_1$ and $\Sigma_yS=Y_2$),
$\sigma:\Sigma_yS\to (\uparrow_a^b)^\perp$ by $\eta\mapsto\uparrow_a^\eta$ is a radial cone-neighborhood isometry,
which implies that $\dim((\uparrow_a^b)^\perp)\geq\dim(\Sigma_yS)$.
On the other hand, by (1.7.1) ``$(\uparrow_a^b)^\perp=(\uparrow_a^b)^{\geq\frac\pi2}$'' implies that $(\uparrow_a^b)^\perp$ is convex in $\Sigma_a(\Sigma_yX)$, so by (1.8.1)
$\dim((\uparrow_a^b)^\perp)\leq\dim(\Sigma_a(\Sigma_yX))-1=\dim(\Sigma_yS).$
Hence, $\dim((\uparrow_a^b)^\perp)=\dim(\Sigma_yS)$, and thus by
(1.10.2) $\sigma$ is surjective, i.e.
$$\Uparrow_a^{\Sigma_yS}=(\uparrow_a^b)^\perp.\eqno{(6.6)}$$
Then by (1.7.1-2), it is not hard to see that (6.5) and (6.6) implies (6.4).

Next, based on the above claim we will complete the proof according to
$[a\Sigma_yS]=(\{a,\tilde
a\}*\widetilde{\Sigma_yS})/\Bbb Z_2$ or $\{a\}*\Sigma_yS$.

\vskip1mm

\noindent Subcase 1: $[a\Sigma_yS]=(\{a,\tilde a\}*\widetilde{\Sigma_yS})/\Bbb
Z_2$ with $\widetilde{\Sigma_yS}/\Bbb Z_2=\Sigma_yS$.

\vskip1mm

Note that it is not hard to
see that $(\{a,\tilde a\}*\widetilde{\Sigma_yS})/\Bbb Z_2$ has an empty boundary.
Then by (6.4), we can apply (6.2.3) to conclude that $\{\xi\}*[a\Sigma_yS]$
can be embedded isometrically into $\Sigma_yX$ (with full dimension)
or $\Sigma_yX=(\{\xi,\widetilde{\xi}\}*\widetilde{[a\Sigma_yS]})/\Bbb Z_2$ with
$\widetilde{[a\Sigma_yS]}/\Bbb Z_2=[a\Sigma_yS]$. In any case,
it has to hold that $\xi=b$, so $\Sigma_yX=(\{\xi,\widetilde{\xi}\}*\widetilde{[a\Sigma_yS]})/\Bbb
Z_2$ (otherwise, plus (6.2-3), $\Sigma_yX=\{\xi\}*[a\Sigma_yS]$ which contradicts
$\partial(\Sigma_yX)=\emptyset$). Since there is a unique minimal geodesic from $a$ to $\xi$,
$a$ has to correspond to a fixed point with respect to the $\Bbb Z_2$-action on $\widetilde{[a\Sigma_yS]}$.
Then it can be seen that, in $\Sigma_yX=(\{\xi,\widetilde{\xi}\}*\widetilde{[a\Sigma_yS]})/\Bbb
Z_2$, $(a)^\perp=[\xi\Sigma_yS]$ and $\Sigma_yX$ is also isometric to $(\{a,\tilde a\}*\widetilde{(a)^\perp})/\Bbb
Z_2$ with $\widetilde{(a)^\perp}/\Bbb Z_2=(a)^\perp$, which implies (6.1.2).

\vskip1mm

\noindent Subcase 2: $[a\Sigma_yS]=\{a\}*\Sigma_yS$.

\vskip1mm

In this case, (6.6) implies that $(\uparrow_a^b)^\perp=\Sigma_yS$ in
$\Sigma_a(\Sigma_yX)$. By (6.2.3) and (6.5),
$\Sigma_a(\Sigma_yX)=\{\uparrow_a^b\}*\Sigma_yS$ or
$(\{\uparrow_a^b,\widetilde{\uparrow_a^b}\}*\widetilde{\Sigma_yS})/\Bbb
Z_2$ with $\widetilde{\Sigma_yS}/\Bbb Z_2=\Sigma_yS$.
In fact, the former case does not occur because $\partial\Sigma_a(\Sigma_yX)=\emptyset$,
and thus $[ab]*\Sigma_yS$ cannot be embedded isometrically into
$\Sigma_yX$ (otherwise $\Sigma_a(\Sigma_yX)=\{\uparrow_a^b\}*\Sigma_yS$). Then
by (6.2.1), ``$\Sigma_a(\Sigma_yX)=(\{\uparrow_a^b,\widetilde{\uparrow_a^b}\}*\widetilde{\Sigma_yS})/\Bbb
Z_2$ and $[a\Sigma_yS]=\{a\}*\Sigma_yS$'' implies that

\vskip1mm

{\it There is an open and dense subset $B$ in
$\Sigma_yS$ such that there are just two minimal geodesics between
an interior point of $[ab]$ and any point in $B$, and any
$[\uparrow_a^b\uparrow_a^\eta]\subset\Sigma_a(\Sigma_yX)$ with
$\eta\in\Sigma_yS$ is realized by a convex spherical surface
bounded by $[a\xi]$, some $[a\eta]$ and $[\xi\eta]$.} \hfill(6.7)

\vskip1mm

\noindent By (6.7), any $\zeta\in\Sigma_yX$ with
$|a\zeta|\leq\frac\pi2$ lies in such a surface because
$\Sigma_a(\Sigma_yX)=[\uparrow_a^b\Sigma_yS]$. This implies the former
part of (6.1.2), i.e. $\Uparrow_{a}^{(a)^\perp}=\Sigma_{a}(\Sigma_yX)$, and that
$$(a)^\perp\subseteq[\xi\Sigma_yS].\eqno{(6.8)}$$
Then we can define a multi-valued map
$$\sigma:\Sigma_yS\to \Sigma_{\xi}(\Sigma_yX) \text{ by }
\eta\mapsto\Uparrow_{\xi}^{\eta}\cap\Sigma_\xi(a)^\perp.$$ Subclaim: {\it
$\sigma$ is a radial cone-neighborhood isometry}.
This together with (6.7) implies that there is a dense subset $B'$ in $\Sigma_yS$ with $B'\subseteq B$ such that
$|\sigma(\eta)|=2$ for all $\eta\in B'$ and $|\sigma(\eta')|=1$
for all $\eta'\in\Sigma_yS\setminus B'$. And the subclaim and (1.7.3) implies that
for any $\zeta,\zeta'\in(a)^\perp$ with $\zeta'$ close to $\zeta$ there is a
$[\zeta\zeta']\subset(a)^\perp$.
We now prove that any $[a\zeta]$ is perpendicular to $[\zeta\zeta']$ (i.e., the latter part of (6.1.2))
according to $|\sigma(\eta)|=2$ or 1 for $\eta\in B'$.

If $|\sigma(\eta)|=2$ for $\eta\in B'$,
it is easy to see that
$$\text{there is a unique minimal geodesic from $a$ to any $\zeta\in(a)^\perp$}.\eqno{(6.9)}$$
Otherwise, there are at least two minimal geodesics from $a$ to some $\zeta\in[\xi\eta]^\circ$ with $\eta\in B'\subseteq B$,
so by (1.7.3) there are at least three convex spherical surfaces containing
$[a\xi]$ and $[a\eta]$, i.e. there are at least three minimal geodesics from
$\uparrow_a^b$ to $\uparrow_{a}^\eta$, which contradicts $\Sigma_a(\Sigma_yX)=(\{\uparrow_a^b,\widetilde{\uparrow_a^b}\}*\widetilde{\Sigma_yS})/\Bbb
Z_2$. Note that because $[\zeta\zeta']\subset(a)^\perp$, by (6.9) and (1.7.3) we can see that
$[a\zeta]$ is perpendicular to $[\zeta\zeta']$.

If $|\sigma(\eta)|=1$ for $\eta\in B'$, then
$\sigma$ is an isometry, i.e. $(a)^\perp=\{\xi\}*\Sigma_yS$, and thus $(a)^\perp$
is convex in $\Sigma_yX$. Then it is not hard
to see that any $[a\zeta]$ is perpendicular to $(a)^\perp$ for $\zeta\in ((a)^\perp)^\circ=(a)^\perp\setminus\Sigma_yS$
by Lemma 1.3 and for $\zeta\in \Sigma_yS\subset(a)^\perp)$ by (1.7.3) (note that there is a
unique minimal geodesic from $a$ to $\zeta\in\Sigma_yS$ because $[a\Sigma_yS]=\{a\}*\Sigma_yS$). It follows that $[a\zeta]$ is perpendicular to $[\zeta\zeta']$.

In the rest, we need only to verify the above subclaim.
If $\xi\neq b$, it is easy to see that $\Uparrow_\xi^{\Sigma_yS}\subseteq(\Sigma_\xi[ab])^\perp$ (cf. (6.2.2)). Then by (1.7.2), $[\xi\Sigma_yS]\subseteq(a)^\perp$, so by (6.8)
$$(a)^\perp=[\xi\Sigma_yS].$$
It follows that $\Uparrow_{\xi}^{\eta}\cap\Sigma_\xi(a)^\perp=\Uparrow_{\xi}^{\eta}$ for any
$\eta\in\Sigma_yS$,
which implies the Subclaim by (6.2.2). If $\xi=b$ (it possibly occurs that $(a)^\perp\subsetneq[\xi\Sigma_yS]$),
the Subclaim holds because  $[\nu\Sigma_yS]$ with $\nu\in[ab]^\circ$ converges to $(a)^\perp$ as $\nu$ converges to $\xi$
and by (6.2.2) $\Sigma_yS\to (\Sigma_{\nu}[ab])^\perp$
with by $\eta\mapsto\Uparrow_{\nu}^{\eta}$ is a radial cone-neighborhood
isometry.
\hfill$\qed$
\enddemo

\demo{Proof of (6.1.3)}

By the structure of $\Sigma_yX$ got in each case in the proof of (6.1.2), one can check that either $(\uparrow_y^x)^\perp=\bigcup_{i=1}^2[\xi_i\Sigma_yS]$
($\xi_1$ may be equal to $\xi_2$) or $(\uparrow_y^x)^\perp\subseteq[\xi\Sigma_yS]$; moreover, $\sigma:\Sigma_yS\to \Sigma_{\xi_i}(\Sigma_yX)$ defined
by $\eta\mapsto\Uparrow_{\xi_i}^{\eta}$ or by $\eta\mapsto\Uparrow_{\xi}^{\eta}\cap\Sigma_\xi(\uparrow_y^x)^\perp$
respectively is a radial cone-neighborhood
isometry  (cf. (6.2.2)). On the other hand, because $\gamma_{\xi_i,y}$ (here $\xi_i$ may be $\xi$)  belongs to some $F_{v_i}$,
by (4.1.3) the map $\sigma_{\xi_i}:\Sigma_yS\to \Sigma_{\xi_i}(\Sigma_yX)$
in (4.3) determined by $F_{v_i}$  is also a radial cone-neighborhood
isometry. By Lemma 5.1, $\sigma=\sigma_{\xi_i}$,
which implies that $(\uparrow_y^x)^\perp=\bigcup_{i=1}^2[\xi_i\Sigma_yS]_{F_{v_i}}\subseteq\Sigma_y\Cal F$.
\hfill$\qed$
\enddemo

%%%%%%%%%%%%%%%%%%%%%%%%%%%%%%%% section 8 %%%%%%%%%%%%%%%%%%%%%

\vskip4mm

\head 7. Reduction on Simply Connectedness and Local Orientability
\endhead

\vskip4mm

Consider the gradient flow of $d(\cdot,\Cal F)$ (cf. [Pet1]). Our goal is to
construct a flow that coincides with the gradient flow on $X\setminus \Cal F$,
and that flows out some $F_v\setminus S$ (see (0.6.2)); which
requires that $F_v\setminus S$ has a neighborhood
$U$ such that $U\setminus F_v$ has two connected components.
We will prove this property by assuming that $X$ is simply connected,
and locally orientable in the sense of [HS].

The purpose of this section is to show that Conjecture 0.3 holds
if it holds on simply connected and locally simply connected spaces.
Hence, without loss of generality we may always assume that $X$ is
simply connected and locally orientable. Note that if $X$ is topologically
nice (cf. [Ka]), then $X$ is always locally orientable, but
the converse does not hold.

For $A\in \op{Alex}^n(\kappa)$, $A$ is said to be {\it locally orientable} ([HS],
[Pet2]) if
$$H^n(A,A\setminus\{a\};\Bbb Z)\cong\Bbb Z, \quad \forall\,\, a\in A.$$
If $A$ is compact and $\partial A=\emptyset$,
then $A$ is said to be {\it orientable} ([HS], [Pet2]) if $A$ is
locally orientable and
$$H^{n}(A,A\setminus\{a\};\Bbb Z)\cong
H^n(A;\Bbb Z), \quad \forall\,\, a\in A.$$
If $A$ is not compact, we will replace the cohomology group by one
with a compact support.

Note that contrary to a topological manifold, that $X$ is simply connected
does not imply that $X$ is orientable ([HS]). According to [HS], if $A$ has an
empty boundary, then $A$ is either orientable or has a (ramified) double
which is orientable ([HS]), i.e., there exists an
orientable Alexandrov space $\tilde A$ with the same dimension and
lower curvature bound which admits an isometric involution $i$ such
that $\tilde A/i$ is isometric to $A$.

\proclaim{Lemma 7.1} {\rm (7.1.1)} Conjecture 0.3 holds if it holds when $X$ is simply connected.

\noindent {\rm (7.1.2)} For the case that $\dim(S)=\dim(X)-2$, Conjecture 0.3 holds
if it holds when $X$ is locally orientable.
\endproclaim

\demo{Proof} (7.1.1) This follows from [Li], where the argument goes through without
the restriction that $X$ is topological nice.

(7.1.2) If $X$ is not locally orientable, then we consider the branched
orientable double cover of $X$, denoted by $\tilde X$. Let $\phi$ be
the covering map from $\tilde X$ to $X$.
Claim: {\it There is a soul $\tilde S$ in $\tilde X$
and a Sharafutdinov retraction $\tilde \pi:\tilde X\to \tilde S$ such that
if $\tilde \pi$ is a submetry then $\pi:X\to S$ is a submetry}.

We now verify the above claim. Let
$\tilde f\triangleq f\circ\phi:\tilde X\to \Bbb R,$
where $f$ is the Busemann function defined in Section 1. The key
is to prove that $\tilde f$ is also a concave function.

Subclaim 1: $\tilde\Omega_c\triangleq\tilde f^{-1}([c,
c_0])\ (=\phi^{-1}(\Omega_c))$ is totally convex in $\tilde X$, where
$c_0=\max f$, and $\partial\tilde\Omega_c=\phi^{-1}(\partial\Omega_c)$.
Then for any $\tilde x\in\tilde\Omega_c$
$$\tilde f(\tilde x)=f(\phi(\tilde x))=|\phi(\tilde
x)\partial\Omega_c|+c=|\tilde x\partial\tilde\Omega_c|+c.$$
It follows that $\tilde f$ is a concave function, so by
the gradient flow of $\tilde f$ (resp. $f$), one get a retraction
$\tilde \pi_0:\tilde X\to\tilde C_0$ (resp. $\pi_0: X\to C_0$),
where $\tilde C_0=\tilde f^{-1}(c_0)$ and $C_0=f^{-1}(c_0)$.

We will prove the above claim according to $\partial \tilde C_0=\emptyset$ or not.

If $\partial \tilde C_0=\emptyset$, we let $\tilde S=\tilde C_0$. Since $\tilde f\triangleq f\circ\phi$,
$\phi$ maps a gradient curve of $\tilde f$ to a gradient curve of $f$,
so $\phi\circ\tilde \pi_0=\pi_0\circ\phi$.
Hence, if $\tilde \pi$ is a submetry, then $\pi_0$ is a submetry.
If $C_0=S$, i.e. $\pi_0=\pi$, then the proof is done. If $C_0\neq S$, then $\pi=\pi_1\circ\pi_0$,
where $\pi_1:C_0\to S$ is the retraction determined by
the gradient flow of $f_1\triangleq d(\cdot, \partial C_0)$ on $C_0$ (note that $S=f_1^{-1}(\max f_1)$, see part b of Section 1).
Note that $\dim(S)=\dim(C_0)-1$ because $\dim(S)=\dim(X)-2$ and $C_0\neq S$. It follows that $\pi_1$ is a submetry ([SY]), so is $\pi$.

If $\partial \tilde C_0\neq\emptyset$, then similar to Subclaim 1 we make the following claim.

\noindent Subclaim 2: If $\partial \tilde C_0\neq\emptyset$,
then $\partial \tilde C_0=\phi^{-1}(\partial C_0)$ and thus $\partial C_0\neq\emptyset$.
Let $\tilde f_1=d(\cdot, \partial\tilde C_0)$ on $\tilde C_0$. Then $\tilde f_1=f_1\circ\phi$ \
(i.e., for any $\tilde x\in\tilde C_0$, $|\tilde x\tilde C_0|=|\phi(\tilde x)C_0|$).
It follows that $\phi$ maps a gradient curve of $\tilde f_1$ to
a gradient curve of $f_1$. Let $\tilde S=\tilde f_1^{-1}(\max \tilde f_1)$.
Note that $\tilde S=\phi^{-1}(S)$, and $\partial\tilde S=\emptyset$ because $\partial S=\emptyset$.
By patching together the gradient flows of $\tilde f$ and $\tilde f_1$,
one get a Sharafutdinov retraction $\tilde \pi:\tilde X\to\tilde S$.
Similarly, $\phi\circ\tilde \pi=\pi\circ\phi$ (note that $\pi=\pi_1\circ\pi_0$),
so if $\tilde \pi$ is a submetry then $\pi$ is a submetry.

In the rest of the proof, we need only to verify Subclaim 1 and 2.

Let  $\sigma(t)|_{[0,+\infty)}$ be a ray in
$X$ with $\sigma^+(0)\in \Uparrow_{\sigma(0)}^{\partial\Omega_{c}}$
and $c<f(\sigma(0))$ (see Lemma 1.2). Observe that, at any $\tilde
x\in\phi^{-1}(\sigma(0))$, there is a ray
$\tilde\sigma(t)|_{[0,+\infty)}$ such that
$\tilde\sigma(t)\in\phi^{-1}(\sigma(t))$ and $\phi|_{\tilde\sigma(t)|_{[0,+\infty)}}$
is an isometry. Note that
$$|\sigma(t_1)\sigma(t_2)|=
|\Omega_{f(\sigma(t_1))}\partial\Omega_{f(\sigma(t_2))}| \text{ for
all } 0\leq t_1<t_2,\eqno{(7.1)}$$ and thus it is easy to see that
$$|\tilde\sigma(t_1)\tilde\sigma(t_2)|=
|\phi^{-1}(\Omega_{f(\sigma(t_1))})\phi^{-1}(\partial\Omega_{f(\sigma(t_2))})|.\eqno{(7.2)}$$

Based on (7.1) and (7.2), we first prove that
$\tilde\Omega_c$ is totally convex in $\tilde X$ (see Subclaim 1).
If this is not true, there
will be $[\tilde p\tilde q]$ with $\tilde p,\tilde q\in\tilde
\Omega_c$ such that $[\tilde p\tilde q]\not\subset\tilde\Omega_c$.
Then we can select $\tilde y\in[\tilde p\tilde q]$ such that $\tilde
f|_{[\tilde p\tilde q]}$ achieves minimum $\bar c<c$ at $\tilde y$
and for any neighborhood $U$ of $\tilde y$ there is $\tilde
y'\in[\tilde p\tilde q]\cap U$ with $\tilde f(\tilde y')>\bar c$.
Let $\sigma(t)|_{[0,+\infty)}$ with $\sigma(0)=\phi(\tilde y)$
and $\tilde\sigma(t)|_{[0,+\infty)}$ with
$\tilde\sigma(0)=\tilde y$ be the rays selected as above. By (7.2),
it has to hold that $|\tilde\sigma^+(0)\uparrow_{\tilde y}^{\tilde
p}| =|\tilde\sigma^+(0)\uparrow_{\tilde y}^{\tilde q}|=\frac\pi2,$
and thus
$$|\phi_*(\tilde\sigma^+(0))\phi_*(\uparrow_{\tilde y}^{\tilde p})|=
|\sigma^+(0)\phi_*(\uparrow_{\tilde y}^{\tilde p})|\leq\frac\pi2
\text{ and } |\sigma^+(0)\phi_*(\uparrow_{\tilde y}^{\tilde
q})|\leq\frac\pi2,\eqno{(7.3)}$$ where $\phi_*$ is the branched
double covering map from $\Sigma_{\tilde y}\tilde X$ to
$\Sigma_{\phi(\tilde y)}X$ induced by $\phi$. On the other hand, by (7.1) it
is not hard to see that
$$|\sigma^+(0)\eta|\geq\frac\pi2\ \forall\ \eta\in\Sigma_{\phi(\tilde y)}\Omega_{\bar c}
\text{ and } |\sigma^+(0)\eta|>\frac\pi2\ \forall\
\eta\in(\Sigma_{\phi(\tilde y)}\Omega_{\bar c})^\circ.
\eqno{(7.4)}$$ Note that $\phi([\tilde p\tilde q])$  is a piecewise minimal
geodesic in $\Omega_{\bar c}$ which is totally
convex in $X$. It follows that $\phi_*(\uparrow_{\tilde
y}^{\tilde p}), \phi_*(\uparrow_{\tilde y}^{\tilde
q})\in\Sigma_{\phi(\tilde y)}\Omega_{\bar c}$, and we can assume that
$\phi_*(\uparrow_{\tilde y}^{\tilde p})$ or $\phi_*(\uparrow_{\tilde
y}^{\tilde q})$ is equal to $\uparrow_{\phi(\tilde y)}^{\phi(\tilde y')}\in(\Sigma_{\phi(\tilde
y)}\Omega_{\bar c})^\circ$. Then (7.4) contradicts (7.3).

Next we prove that $\partial\tilde\Omega_c=\phi^{-1}(\partial\Omega_c)$
(so Subclaim 1 is verified). In fact, it is clear that
$\partial\tilde\Omega_c\subseteq\phi^{-1}(\partial\Omega_c)$;
and for $x\in\partial\Omega_c$, by (7.1) and (7.2) we
can conclude that any $\tilde x\in\phi^{-1}(x)$ is not an
interior point of $\tilde\Omega_c$, i.e. $\partial\tilde\Omega_c\supseteq\phi^{-1}(\partial\Omega_c)$.

As for Subclaim 2, we will omit its proof because the proof
is similar to that for ``$\partial\tilde\Omega_c=\phi^{-1}(\partial\Omega_c)$''
(here, the ray $\sigma$ (resp. $\tilde \sigma$) in (7.1) and (7.2)
will be minimal geodesic from $x\in C_0$ to $\partial C_0$ (resp.
from $\tilde x\in\phi^{-1}(x)$ to $\partial\tilde C_0$)). \hfill$\qed$
\enddemo

\remark{Remark {\rm 7.2}} (7.2.1) We point it out that from the above proof, (7.1.2)
still holds (without dimension restriction) for the case that $\partial \tilde C_0
=\emptyset$ and $\dim(S)=\dim(\tilde C_0)-1$. Furthermore, we point it out that
the proof of (7.1.2) implies (7.1.1).

(7.2.2) Consider the case of Conjecture 0.3 for $\dim(S)=\dim(X)-2$:
we may assume that $X$ is locally orientable and simply conmnected; by
passing to a branched cover, we may assume $X$ is locally orientable,
and further pass to the universal cover of $X$ (which has the same
local topology).
\endremark

%%%%%%%%%%%%%%%%%%%%%%%%%%%%%%%% section 8    %%%%%%%%%%%%%%%%%%%%%

\vskip4mm

\head 8. Structures of $S$ and $\tilde S_c$ When
$\dim(X)=4$
\endhead

\vskip4mm

In the rest of the paper, we will be confined to the case that $\dim(X)=4$ and $\dim(S)=2$.
Consider the Sharafutdinov retraction, $\pi|_{\tilde S_c}: \tilde S_c\to S$, where $\tilde
S_c=F_v\cap\partial\Omega_c$ with $c<c_0$, and the inverse map of $\pi|_{\tilde S_c}$, $\varphi_c: S\to \tilde S_c$,
is a radial cone-neighborhood isometry (see Section 4) such that for $q\in S$, $|\varphi(q)|\le k$ and ``$=$'' iff $q\in S_k(F_v)$.

Our goal in this section is to establish the following structural results
on $S$ and $\tilde S_c$, where we always assume the assumptions in Theorem A
and $\dim(S)=2$.

\proclaim{Theorem 8.1} The following holds:

\noindent{\rm (8.1.1)} $S\setminus S_k(F_v)$ is a finite set.

\noindent{\rm (8.1.2)} $\tilde S_c\in \op{Alex}^2(0)$ with
$\partial\tilde S_c=\emptyset$.

\noindent{\rm (8.1.3)}  If $S$ is simply connected, then $\tilde
S_c$ is simply connected.
\endproclaim

The following corollary will be used in the proofs of (0.6.2) and Theorem 0.7.

\proclaim{Corollary 8.2} Let $Q=S\setminus S_k(F_v)$. Then

\noindent{\rm  (8.2.1)} $\tilde Q_c\triangleq\varphi_c(Q)$ satisfies that $|\tilde Q_c|<k\cdot |Q|<\infty$.

\noindent{\rm  (8.2.2)} $\pi|_{\tilde S_c\setminus\tilde Q_c}:
\tilde S_c\setminus \tilde Q_c\to S\setminus Q$ is a metric $k$-cover.

\noindent{\rm (8.2.3)} $\pi|_{\tilde S_c}$ is a local radial isometry at points in $\tilde Q_c$;
and for any $x\in \tilde S_c$ and $r$, we have that
$\pi|_{\tilde S_c}(B_r(x))=B_r(\pi(x))$.
\endproclaim

\demo{Proof} (8.2.1) follows from (8.1.1) and that $\varphi_c$ is a radial cone-neighborhood isometry (see (4.1.1)), and the latter also implies (8.2.2) and (8.2.3).
\hfill$\qed$
\enddemo

We need some preparation for the proof of Theorem 8.1.
Note that for $q\in S$, $(\Sigma_qS)^\perp\neq\emptyset$ (see (1.1)), and $\Sigma_q S$ is a convex circle in $\Sigma_qX$ because $\dim(S)=2$
and $S$ is convex in $X$.

\proclaim{Lemma 8.3} For $q\in S$ and $\eta\in
(\Sigma_qS)^\perp$, $\Uparrow_{\eta}^{\Sigma_qS}$ is a locally convex circle with
perimeter $\leq 2\pi$ in $\Sigma_qX$, and the map $\Uparrow_{\eta}^{\Sigma_qS}\to
\Sigma_qS$  by $\Uparrow_{\eta}^\zeta\mapsto\zeta$ is a
metric cover.
\endproclaim

Lemma 8.3 has the following corollary.

\proclaim{Corollary 8.4} Let $q\in S$. If $\Sigma_qS$ has diameter bigger than
$\frac\pi2$, then

\noindent{\rm (8.4.1)} $[(\Sigma_qS)^\perp\Sigma_qS]=(\Sigma_qS)^\perp*\Sigma_qS$.

\noindent{\rm (8.4.2)} $q$ belongs to $ S_k(F_v)$.
\endproclaim

\demo{Proof} Since $\op{diam}(\Sigma_qS)>\frac\pi2$, for any $\eta\in(\Sigma_qS)^\perp$ the metric cover
$\Uparrow_{\eta}^{\Sigma_qS}\to \Sigma_qS$
in Lemma 8.3 has to be an isometry. This implies that there is a unique
minimal geodesic from $\eta$ to any $\xi\in\Sigma_pS$, and thus (8.4.1) holds.

Let $q'\in S$ be another point with $\op{diam}(\Sigma_{q'}S)>\frac\pi2$.
By (8.4.1), we can use the same arguments as proving (3.1.1) to conclude that:
{\it For any $[qq']$, $\phi_{[qq']}$ is an isometry} (see (3.1.1)). Then by the construction of $F_v$,
we can see (8.4.2).
\hfill$\qed$
\enddemo

Since $\dim(\Sigma_qX)=3$, Lemma 8.3 is
an immediate corollary of Lemma 8.5 below.

\proclaim{Lemma 8.5} Let $Y\in \text{\rm Alex}^3(1)$, and $A$ be a
locally convex circle in $Y$, and $p$ be $\frac\pi2$-apart from $A$
in $Y$. Then $\Uparrow_{p}^{A}$ is a locally convex circle in
$\Sigma_pY$ with perimeter $\leq2\pi$, and the map
$\Uparrow_{p}^{A}\to A$ by $\Uparrow_{p}^a\mapsto a$ is a metric cover.
\endproclaim

By Lemma 1.3 and (1.8.1), we know that
$A^{\geq\frac\pi2}=A^{=\frac\pi2}$ and $\dim(A^{=\frac\pi2})\leq 1$.
Note that if $\dim(A^{=\frac\pi2})=1$ and
$p\in(A^{=\frac\pi2})^\circ$, then Lemma 8.5 is a corollary of
(6.2.2).

\demo{Proof} Since $|pa|=\frac\pi2$ for all $a\in A$ (by
Lemma 1.3), given $[a_1a_2]\subset A$ with $a\in[a_1a_2]^\circ$, any $[pa]$ belongs to a
convex spherical surface spanned by $p$ and $[a_1a_2]$ (by (1.7.3)).
Note that $[pa]$ belongs to a unique such surface (otherwise, there is $[pa_1]$ and
$[pa_2]$ such that $|\uparrow_p^{a_1}\uparrow_p^{a_2}|<|a_1a_2|$, which is impossible
because $|\uparrow_p^{a_1}\uparrow_p^{a_2}|\geq|a_1a_2|$ by (1.7.2)).
Plus $\Sigma_pY\in\text{Alex}^2(1)$, it follows that each
component of $\Uparrow_{p}^{A}$ either is a locally convex circle,
or converges to two locally convex circles which do not intersect
each other. This together with Lemma 5.3 implies that $\Uparrow_{p}^{A}$
is just a locally convex circle, and thus the map $\Uparrow_{p}^{A}\to A$ by
$\Uparrow_{p}^a\mapsto a$ is a metric cover. And by Lemma 8.6 below
$\Uparrow_{p}^{A}$ has perimeter $\leq 2\pi$.
\hfill$\qed$
\enddemo

\proclaim{Lemma 8.6} Let $Y\in \op{Alex}^2(1)$, and let $S^1$ be a locally convex circle in $Y$.
Then $S^1$ has perimeter $\leq 2\pi$.
\endproclaim

\demo{Proof}  Note that $Y$ is homeomorphic to $\Bbb S^2$ if $Y$ is
simply connected ([BGP]), and thus $S^1$ divides $Y$
into two spaces $\in\text{Alex}^2(1)$ containing
$S^1$ as boundary. It then is not hard to see that $S^1$ is
of length $\leq 2\pi$ (essentially by (1.7.2)). If $Y$ is not simply connected, we can
consider its universal covering space to draw the conclusion.
\hfill$\qed$
\enddemo

We are now ready to prove Theorem 8.1.

\demo{Proof of (8.1.1)}

According to (8.4.2), it suffices to show that there
is only a finite number of points in $S$ whose spaces of directions
have diameter $\leq\frac\pi2$. Note that this follows
from that $S$ is compact and a basic fact on an Alexandrov space $Z$
with lower curvature bound ([BGP]): {\it Given any $z\in Z$, for any
$\epsilon>0$ there is a neighborhood $U_\epsilon$ of $z$ such that
the diameter of $\Sigma_{z'}Z$ is bigger than $\pi-\epsilon$ for any
$z'\in U_\epsilon\setminus\{z\}$} (we can prove this by a standard
limiting argument). \hfill$\qed$
\enddemo

\demo{Proof of (8.1.2)}

By (8.1.1), we can let
$S\setminus S_k(F_v)=\{q_1,\cdots,q_m\}$. For convenience, we
let $q$ be an arbitrary point in $\{q_1,\cdots,q_m\}$, and let
$\tilde q$ be an arbitrary point in $\varphi_c(q)$. Since
$\varphi_c$ is a radial cone-neighborhood isometry ((4.1.1)),
$\pi|_{\tilde S_c}$ is a metric $k$-cover on
$\varphi_c( S_k(F_v))$ which implies that $\varphi_c( S_k(F_v))$ is an open
$2$-dimensional Alexandrov space with curvature
$\geq0$; and  $\pi|_{\tilde S_c}$ is
a branched and radial isometric $g$-cover around $\tilde q$
with $\tilde q$ being branched point and $g\leq k$.

It remains to show that $\tilde S_c$ is also of curvature $\geq0$ at
$\tilde q$ in the Alexandrov sense and that $\tilde q$ is also an
interior point. For small $\epsilon>0$, by (1.9.4)
$B_\epsilon(q)\subset S$ can be divided into small
sectors $\{\frak{S}_j\}_{j=1}^l$ (note that $\dim(S)=2$)  with each $\frak{S}_j\cap\frak{S}_{j+1}$
(where $\frak{S}_{l+1}=\frak{S}_1$) being a minimal geodesic of length
$\epsilon$ starting from $q$ such that $B_{\epsilon}(\tilde q)\subset\tilde
S_c$ can be divided into small
sectors $\{\widetilde{\frak S_j}\}_{j=1}^{gl}$,
where $\widetilde{\frak S_{il+l'}}$ is isometric to $\frak S_{l'}$ for
any $0\leq i\leq g-1$ and $1\leq l'\leq l$ and each
$\widetilde{\frak{S}_j}\cap\widetilde{\frak{S}_{j+1}}$ (where
$\widetilde{\frak{S}_{gl+1}}=\widetilde{\frak{S}_1}$) is a minimal
geodesic of length $\epsilon$ starting from $\tilde q$. On the other
hand, $\Sigma_{\tilde q}\tilde S_c=\text{\rm
D}\varphi_c(\Sigma_qS)$ (see (4.1.2)), where  $\text{\rm
D}\varphi_c$ is a radial cone-neighborhood isometry from $\Sigma_qS$
to $\Sigma_{\tilde q}\partial\Omega_c$ of dimension 2.
Hence, $\Sigma_{\tilde q}\tilde S_c$
is a union of several disjoint locally convex circles (see the beginning of the proof of Lemma 4.3).
It therefore follows from Lemma 5.3 and 8.6 $\Sigma_{\tilde q}\tilde S_c$
is a circle of length $\leq 2\pi$. This implies that
$\tilde S_c$ is of curvature $\geq0$ at $\tilde q$, and $\tilde q$ is an interior
point.  (Note that this process implies that $\Bbb Z_g$ can act on
$B_{\epsilon}(\tilde q)$ by isometries such that
$B_{\epsilon}(q)=B_{\epsilon}(\tilde q)/\Bbb
Z_g$.)\hfill$\qed$
\enddemo

Note that the above proof relies on that $\dim(S)=2$. In higher
dimension, we conjecture that $\tilde S_c$
is also an Alexandrov space with curvature $\geq0$ and an empty
boundary.

\remark{Remark \rm 8.7} From the proof of (8.1.2) and Corollary 4.2, one may
conclude that the natural map $\tilde S_c\to \tilde
S_{c'}$ by $\gamma_{w,q}\cap\tilde S_c\mapsto \gamma_{w,q}\cap\tilde S_{c'}$
is an isometry, where $c,c'<c_0$ and $\gamma_{w,p}\subset F_v$ with $p\in S$ and $w\in\Uparrow_p^{\partial\Omega_c}$.
That is, $F_v\setminus S$ with the intrinsic
metric is isometric to $\tilde S_c\times(0,+\infty)$
(however, with the extrinsic metric $F_v\setminus S$ may
not be isometric to $\tilde S_c\times(0,+\infty)$).
\endremark

\demo{Proof of (8.1.3)}

By (8.1.2), $\tilde S_c\in\op{Alex}^2(0)$ and $\partial \tilde S_c=\emptyset$,
so $\tilde S_c$ is a closed surface (cf. [BGP]). Hence, if $\tilde S_c$ is not simply connected,
then $\tilde S_c$ contains a loop $\sigma$ which is not homotopy to a point such that
$\tilde S_c\setminus\sigma$ is still connected and thus $S\setminus
\pi|_{\tilde S_c}(\sigma)$ is connected too. On the other hand, note that we
can assume that $\sigma$ lies in $\varphi_c( S_k(F_v))$ because
$|\varphi_c(S\setminus S_k(F_v))|<\infty$ (see (8.2.1)),
so $\pi|_{\tilde S_c}(\sigma)$ is also a loop in $S$ because $\pi|_{\tilde S_c}$ is a covering
map on $\varphi_c( S_k(F_v))$ (see (8.2.2)). Moreover,
note that $S$ is a
2-dimensional sphere because the simply connected $S$ is also a closed surface,
so it follows that $S\setminus \pi|_{\tilde S_c}(\sigma)$ is not connected, a
contradiction. \hfill$\qed$
\enddemo

%%%%%%%%%%%%%%%%%%%%%%%%%%%%%%%% section 9 %%%%%%%%%%%%%%%%%%%%%

\vskip4mm

\head 9. Proof of (0.6.2)
\endhead

\vskip4mm

We need the following notion: $F_v$ is said to satisfy the Separation Property, if
for any $p\in F_v\setminus S$, there is $B_{r_p}(p)\subset X\setminus S$ such that
$F_v\setminus S$ separates $\bigcup_{p\in F_v\setminus S}B_{r_p}(p)$
into two components, $U^1$ and $U^2$.

\proclaim{Lemma 9.1} Suppose that $F_v$ satisfies the Separating Property with
$B_{r_p}(p)$ and $U^i$ defined in the above. Then the following holds.

\vskip1mm

\noindent{\rm (9.1.1)} Let $U_{p}^i\triangleq U^i\cap B_{r_p}(p)$.
If $U^i_{\bar p}\cap \Cal F=\emptyset$ ($i=1$ or $2$) for some $\bar p\in F_v\setminus S$,
then for any $p\in F_v\setminus S$ and $r_p$ sufficiently small,
$U^i_{p}\cap \Cal F=\emptyset$, i.e., $U^i\cap \Cal F=\emptyset$.

\vskip1mm

\noindent{\rm (9.1.2)} Let $x\in X\setminus\Cal F$, and let $y\in \Cal F$ such that $|xy|=|x\Cal F|$. If $y\in F_v\setminus S$, then any $[xy]$ satisfies that $[xy]\cap U^i_{y}\neq\emptyset$
for $i=1$ or $2$, say $1$
(so $[xy]\cap U^2_{y}=\emptyset$ for small $r_y$); and we can
choose $r_y$ so small that $U^1_{y}\cap \Cal F=\emptyset$, and thus
$U^1\cap \Cal F=\emptyset$.

\vskip1mm

\noindent {\rm (9.1.3)}  Let $x\in X\setminus\Cal F$, and let $y\in \Cal F$ such that $|xy|=|x\Cal F|$. If $y\in F_v\setminus S$, then via gradient curves of $d(\cdot,\Cal F)$ we
can define a distance non-increasing flow $\Psi_t|_{[0,+\infty)}:  F_v\to X$ such that $x=\Psi_{|x\Cal F|}(y)$.
\endproclaim

Note that (0.6.2) is a consequence of the following fact and (9.1.3).

\proclaim{Lemma 9.2} Assume that $X$ is simply connected and locally orientable, $\dim(X)=4$ and
$\dim(S)=2$. Then $F_v$ satisfies the Separating Property.
\endproclaim

In the rest of this section, we will prove Lemma 9.1 and 9.2.

\demo{Proof of Lemma 9.1}

\noindent(9.1.1) Let $\bar p\in F_v\setminus S$ such that $U^1_{\bar p}\cap \Cal F=\emptyset$.
If (9.1.1) is not true, then
there exists $p\in F_v\setminus S$ and $p_j\in \Cal F\cap U^1_{p}$ such that $p_j\to p$ as
$j\to\infty$. Note that $p_j$ lies in some $F_{v_j}$,
where $v_j\in \Uparrow_{p_0}^{\partial\Omega_{c}}$ (for $p_0$ refer to the definition of $F_v$ in Section 3).
Since $p_j\in U_p^1$ and $U_p^1\cap F_v=\emptyset$, it holds that $F_{v_j}\neq F_v$.
On the other hand, we can assume that $F_{v_j}$ converges to
some $F_{w}\ni p$ with $w\in \Uparrow_{p_0}^{\partial\Omega_{c}}$ as $j\to\infty$.
Note that $F_w=F_v$ by (3.1.3) because $p\in(F_w\cap F_v)\setminus S$.
Then we claim that $F_{v_j}=F_v$, so we get a contradiction and the proof is done.
Observe that ``$U^1_{\bar p}\cap \Cal F=\emptyset$'' implies that $U^2_{\bar p}\cap F_{v_j}\neq\emptyset$ for large $j$.
Let $\bar p_j\in U^2_{\bar p}\cap F_{v_j}$. Since $F_{v_j}\to F_v$, for large $j$ there is
$[p_j\bar p_j]_{F_{v_j}\setminus S}$, a shortest path from $p_j$ to $\bar p_j$ in $F_{v_j}\setminus S$, such that
$[p_j\bar p_j]_{F_{v_j}\setminus S}\subset U^1\cup U^2\cup F_v$.
It follows that $[p_j\bar p_j]_{F_{v_j\setminus S}}\cap F_v\neq\emptyset$, so
$F_{v_j}=F_v$ by (3.1.3), i.e. the claim is verified.

\vskip2mm

\noindent(9.1.2) Let $x\in X\setminus\Cal F$, and let $y\in F_v\setminus S$
such that $|xy|=|x\Cal F|$. Note that any $[xy]$ satisfies that $[xy]\setminus\{y\}\cap \Cal F=\emptyset$,
so by the Separating Property  we can assume that $[xy]\cap U_y^1\neq\emptyset$
(and $[xy]\cap U^2_{y}=\emptyset$ for small $r_y$). Then by (9.1.1),
it suffices to show that $U^1_{y}\cap \Cal F=\emptyset$ for small $r_y$.
If it is not true, then there is $p_j\in \Cal F\cap U^1_{y}$ such that $p_j\to y$ as
$j\to\infty$. Similar to the proof of (9.1.1),
$p_j$ lies in some $F_{v_j}$ with $F_{v_j}\neq F_v$ (i.e. $F_{v_j}\cap F_v=S$),
and thus we can assume that $U^2_{y}\cap F_{v_j}=\emptyset$.
Then we claim that $[xy]\cap F_{v_j}\neq\emptyset$ for large $j$,
which contradicts ``$|xy|=|x\Cal F|$'', and thus $U^1_{y}\cap \Cal F=\emptyset$  for small $r_y$.

We now need only to verify the claim right above. We first observe
that the proof of (0.6.1) for $y\in S$ in Section 6
implies that the distance function $d(\cdot, F_{v_j})$ is concave
at any $x'\in X\setminus F_{v_j}$ whose nearest point in $F_{v_j}$ lies in $F_{v_j}\setminus S$.
Note that $y\notin S$ and $y$ is close to $F_{v_j}$, so there is $x_j\in [xy]$
with $|x_jy|\geq |yF_{v_j}|$ such that each point in $[x_jy]$ has its nearest point
in $F_{v_j}$ lie in $F_{v_j}\setminus S$.
Hence, if the claim is not true, i.e. $[xy]\cap F_{v_j}=\emptyset$, then
$d(\cdot, F_{v_j})$ is concave on $[x_jy]$ with $|x_jF_{v_j}|\geq |x_j\Cal F|=|x_jy|\geq |yF_{v_j}|$.
It follows that $$|yF_{v_j}|=\min\{|zF_{v_j}||\ z\in[x_jy]\}.\eqno{(9.1)}$$
On the other hand, from the proof of (0.6.1) for $y\in S$ in Section 6
we know that $(\uparrow_y^x)^\perp=\Sigma_y F_v$
and $B_{\frac\pi2}(\uparrow_y^x)=[\uparrow_y^x\Sigma_y F_v]$. Plus
$[xy]\cap U^1_{y}\neq\emptyset$, we can easily see that any $\uparrow_y^{x'}$ with $x'\in U^1_y$
satisfies that $|\uparrow_y^x\uparrow_y^{x'}|<\frac\pi2$. Note that
the nearest point of $y$ in $F_{v_j}$, $\bar y$, lies in $U^1_y$ because
$p_j\to y$ and $U^2_{y}\cap F_{v_j}=\emptyset$, and thus
$|\uparrow_y^x\uparrow_y^{\bar y}|<\frac\pi2$. Hence, by the first variation formula,
for $z\in[x_jy]\subset[xy]$ close to $y$ enough $$|zF_{v_j}|\leq |z\bar y|<|\bar yy|=|yF_{v_j}|,$$
which contradicts (9.1). So, the above claim follows.

\vskip2mm

\noindent (9.1.3) Let $x\in X\setminus\Cal F$, and let $y\in F_v\setminus S$ such that $|xy|=|x\Cal F|$.
By (9.1.2), given a $[xy]$,
we can assume that  $[xy]\cap U^1\neq\emptyset$ and $U^1\cap \Cal F=\emptyset$.

By (0.6.1), the distance function $d(\Cal F,\cdot)$ is concave on $X\setminus \Cal F$.
Then starting from any $\bar x\in U^1$, there is a $d(\Cal F,\cdot)$-gradient curve ([Pet1]).
Note that for any $p\in F_v\setminus S$, by (9.1.1) we can let $r_p$ be so
small that $|\bar x\Cal F|=|\bar xF_v \setminus S|$ and thus
the $d(\Cal F,\cdot)$-gradient curve at $\bar x$ does not pass through $U^2\cup F_v\setminus S$.
And note that the gradient flow of the concave function $d(\Cal F,\cdot)$ (on $X\in\op{Alex}(0)$)
preserves distance non-increasing ([Pet1]). Then using the limiting argument, at any $p\in F_v$
we can construct a unique curve $\varsigma_p(t)|_{[0,+\infty)}$
($\varsigma_p(0)=p$) such that
$\varsigma_p(t)|_{(0,\epsilon)}\subset U^1$ for some $\epsilon>0$
and $\varsigma_p(t)|_{[t_0,+\infty)}$ with $t_0>0$ is the $d(\Cal
F,\cdot)$-gradient curve starting from $\varsigma_p(t_0)$. Note that
$[xy]=\varsigma_y(t)|_{[0,|x\Cal F|]}$. Hence, we can define a flow
$$\Psi_t|_{[0,+\infty)}:F_v\to X \text{ by } p\mapsto \varsigma_p(t) \eqno{(9.2)}$$
such that $\Psi_{|x\Cal F|}(y)=x$. (Note that for $p\in S$, it may occurs that $\Psi_t(p)=p$.)

We now show that $\Psi_t$ preserves distance non-increasing.
Since the gradient flow of $d(\Cal F,\cdot)$
preserves distance non-increasing, it suffices to show that
any shortest path in $F_v$ is a piecewise minimal geodesic in $X$.
In fact, this follows from $|S\setminus S_k(F_v)|<\infty$ (see (8.1.1)) and Corollary 4.2.
\hfill$\qed$
\enddemo

Let's first observe a consequence of the local orientability in dimension $4$.

\proclaim{Lemma 9.3} Assume that $X$ is locally orientable, $\dim(X)=4$ and
$\dim(S)=2$. Let $p$ be an arbitrary point in $\Cal
F\setminus S$. Then

\noindent{\rm (9.3.1)} $\Sigma_pX$ is homeomorphic to $\Bbb S^3$.

\noindent{\rm(9.3.2)} If $p$ lies in some $F_v$,
then $\Sigma_p F_v$ is homeomorphic to $\Bbb S^2$.
\endproclaim

\demo{Proof}  Let $c=f(p)\ (<c_0)$. Since $p\in\Cal F\setminus S$, $p$ is an interior point
of a ray in $\Cal F$. Then it is easy to see that $\Sigma_pX$ is a
spherical suspension over cross section $\Sigma_p(\partial\Omega_c)\in\text{Alex}^2(1)$ without
boundary (note that $\partial\Sigma_pX=\emptyset$).

(9.3.1) It suffices to show that $\Sigma_p(\partial\Omega_c)$ is a sphere. Since $\Sigma_p(\partial\Omega_c)\in\text{Alex}^2(1)$,
$\Sigma_p(\partial\Omega_c)$ is a closed surface ([BGP]), and thus
$\Sigma_p(\partial\Omega_c)$
is a sphere or projective space because of the positive curvature (ref. Theorem 1.8 in [Ma]).
Note that $H^3(\Sigma_pX;\Bbb Z)\cong\Bbb Z$ because $X$ is locally orientable.
This together with that $\Sigma_pX$ is a spherical suspension over $\Sigma_p(\partial\Omega_c)$
implies that $\Sigma_p(\partial\Omega_c)$ is a sphere.

(9.3.2) Let $\tilde S_c\triangleq F_v\cap\partial\Omega_c\ni p$.
Note that $\Sigma_p F_v$ is a
spherical suspension over $\Sigma_p\tilde S_c$ (see the contents after (4.2)), and
$\Sigma_p\tilde S_c$ is a locally convex circle  in $\Sigma_p(\partial\Omega_c)$
(see the proof of (8.1.2)). It follows that $\Sigma_p F_v$ is
homeomorphic to $\Bbb S^2$. \hfill$\qed$
\enddemo

\demo{Proof of Lemma 9.2}

We will use the following fundamental result on the topology of
$Y\in\op{Alex}(\kappa)$ (ref. 4.4 in [Per2]): {\it Given $y\in
Y$, there is a small $r>0$ such that $B_r(y)\subset Y$ is
homeomorphic to $B_r(O)\subset T_y$, where $T_y$ is the
tangent cone at $y$ with
vertex $O$, and that $\partial B_{r'}(y)$ is homeomorphic to
$\Sigma_yY$ for all $0<r'\leq r$}.

Since $X$ is simply connected, $S$ has to be simply
connected; and thus each $\tilde S_c$ $(c<c_0)$ is also simply connected
by (8.1.3), so is $F_v\setminus S$ (by Remark 8.7). Then it suffices to show that
for small $r_p$, $B_{r_p, F_v}(p)$ separates $B_{r_p}(p)\ (\subset X\setminus S)$ into two components,
where $B_{r_p, F_v}(p)$ denotes the closed ball in $F_v$ with center $p$ and radius $r_p$.
By Corollary 4.2 and 5.2, $r_p$ can be chosen such that $\partial B_{r_p, F_v}(p)\subset\partial B_{r_p}(p)$.
Then applying Lemma 9.3 and the above topological property to $X$ and $F_v\setminus S$ at $p$
(note that $F_v\setminus S$ belongs to Alex$^3(0)$ by (8.1.2) and Remark 8.7), we can conclude that
$B_{r_p, F_v}(p)$ separates $B_{r_p}(p)$ into two components for small $r_p$.
\hfill$\qed$
\enddemo

%%%%%%%%%%%%%%%%%%%%%%%%%%%%%%%% section 10 %%%%%%%%%%%%%%%%%%%%%

\vskip4mm

\head 10. Proof of Theorem 0.7 and (A1)
\endhead

\vskip4mm

In this section, we will first complete the proof of (A1).

\demo{Proof of (A1) by assuming Theorem 0.7}

From (0.4.1) and (0.4.2), we will assume that $\dim(S)=2$. We need to show that, for all $x\in X$ and small $r$, $\pi(B_r(x))=B_r(\pi(x))$. By Lemma 1.5, we only need to consider $x\in
X\setminus\Cal F$. Let $y\in \Cal F$ such that $|xy|=|x\Cal F|$.

If $y\not\in S$,  then $y$ belongs to a unique $F_v$ by (3.1.3). By Lemma 7.1 and (7.2.2),
without loss of generality we may assume that $X$ is simply
connected and locally orientable, and thus $F_v$ satisfies the Separating Property (Lemma 9.2).
Then by (9.1.3),
there is a distance non-increasing flow $\Psi_t|_{[0,+\infty)}:  F_v\to X$ in (0.6.2) such that $x=\Psi_{|x\Cal F|}(y)$.
Then by Theorem 0.7, for all $r$
$$\pi(B_r(x))=B_r(\pi(x)).\eqno{(10.1)}$$

If $y\in S$, by the continuity of $\pi$, (10.1) still holds (i.e. (A1) follows) once we show that
$\{x\in X\setminus\Cal F| \text{ if $y\in\Cal F$ satisfies $|xy|=|x\Cal F|$,
then $y\in S$}\}$ is of zero-measure. In fact, this
follows from the fact: if $y$ lies in $S$, then
by (6.1.1) $|\uparrow_y^x\Uparrow_y^{\partial\Omega_c}|=\frac\pi2$ for any
$[xy]$ and thus $\Uparrow_y^x$ consists of at most two elements.
\hfill$\qed$
\enddemo

In the rest of paper, we will prove Theorem 0.7. Consider the following diagram:
$$\CD @. \tilde S_c \\ @.  @VV \pi V \\
\tilde S_c @> \pi\circ \Psi_t >> S,
\endCD$$
where $\Psi_t$ is defined in (0.6.2) which is distance non-increasing.
Because $\pi|_{\tilde S_c}:\tilde S_c\to S$ is a branched cover,
(see (8.2.2-3)), apriori $\pi\circ\Psi_t$ may not have any lifting map.
However, we find that in our situation (see Lemma 10.1 below),
$\pi\circ\Psi_t$ does have a lifting map $\widetilde{\pi\circ\Psi_t}:\tilde S_c\to
\tilde S_c$. Note that $\pi\circ\Psi_t$ is distance non-increasing
because both $\pi$ and $\Psi_t$ are distance non-increasing, so it is easy to see that
$\widetilde{\pi\circ\Psi_t}$ is also distance non-increasing. Moreover, we can prove that
$\widetilde{\pi\circ\Psi_t}$ is surjective, so it is an isometry (by [Pet2])
which implies Theorem 0.7.

Obviously, around the points in
$S\setminus S_k(F_v)$ (i.e., the branched points of
$\pi|_{\tilde S_c}$), the ``lifting'' work will be hard. In order to
do it, it is important to make clear the structure of
$$\tilde Q_c^t\triangleq\left\{\tilde p\in\tilde S_c|
\pi(\Psi_t(\tilde p))\not\in S_k(F_v)\right\}.$$
Note that $\tilde Q_c^0=\tilde Q_c=\varphi_c(S\setminus S_k(F_v))$, which is a finite set by (8.2.1).

\proclaim{Lemma 10.1} For small $t$, $|\tilde Q_c^t|=|\tilde Q_c|$
and $\tilde Q_c^t$ converges to $\tilde Q_c$ as $t\to 0$.
\endproclaim

In its proof, besides the distance non-increasing property of $\pi$
and $\Psi_t$, the following lemma on 2-dimensional
topology  will play a key role.

\proclaim{Lemma 10.2} Let $\rho: \Bbb S^2\to \Bbb S^2$ be a continuous map of degree
$k\geq1$. Then $\rho$ is onto, and the set
$\{p\in\Bbb S^2|\ |\rho^{-1}(p)|<k\}$ contains at most only finite points.
\endproclaim

Here, that $\rho$ is onto is a standard consequence of a non-vanishing degree,
so Lemma 10.2 is trivial when $k=1$.
However, we do not find in literature on the finiteness of the set of points whose inverse
image is less than $k$ when $k\geq 2$.

\demo{Proof} Let $\rho_*:H_2(\Bbb S^2;\Bbb Z)\to H_2(\Bbb S^2;\Bbb
Z)$ be the homomorphism induced by $\rho$, and assume that
$\rho_*(e)=k\cdot e$, where $e$ is the generator of $H_2(\Bbb
S^2;\Bbb Z)\cong\Bbb Z$. Recall that $k$ is called the degree of
$\rho$, and ``$k\geq 1$'' implies that $\rho$ has to be an onto map.

In order to see that $L\triangleq\{p\in\Bbb S^2|\ |\rho^{-1}(p)|<k\}$
is a finite set, it suffices to show that $L$ consists of isolated points,
i.e. for any $p\in L$ there is a
neighborhood $U$ of $p$ such that $|\rho^{-1}(q)|\geq k$ for any
$q\in U\setminus\{p\}$.

We set $\rho^{-1}(p)=\{\tilde p_1,\cdots,\tilde p_l\}\
(l<k)$, and let $\tilde B_i$ be disjoint (open) discs with center
$\tilde p_i$, and let $B$ be a disc with center $p$ such that
$\rho(\tilde B_i)\subseteq B$. Then we let $A$ (resp. $\tilde A$)
denote $\Bbb S^2\setminus\{p\}$ (resp. $\Bbb S^2\setminus\{\tilde
p_1,\cdots,\tilde p_l\}$). Note that $A\cap B$ and each $\tilde
A\cap\tilde B_i$ are homeomorphic to a cylinder, and
$$\rho(\tilde A)=A \text{ and } \rho(\tilde A\cap\tilde B_i)\subseteq A\cap B.$$
Then by considering the following commutative Mayer-Vietoris sequences
$$\CD
0 @>>> H_2(\Bbb S^2; \Bbb Z)& @>>> &H_1(\tilde A\cap \bigcup_{i=1}^l\tilde B_i; \Bbb Z) & @>>> &
H_1(\tilde A; \Bbb Z) & @>>> 0  \\
   & &@V\rho_* VV &                 & @V \rho_* VV &  & @ V \rho_* VV&  \\
0 @>>> H_2(\Bbb S^2; \Bbb Z)& @>>> &H_1(A\cap B; \Bbb Z) & @>>> &0&  @>>>0 &
\endCD$$
we can conclude that there exists $k_i$ such that
$$\rho_*(\sum_{i=1}^l\tilde e_i)=\sum_{i=1}^l\rho_*(\tilde e_i)=\sum_{i=1}^l(k_i\cdot \bar e)=k\cdot \bar e, \eqno{(10.2)}$$
where $\tilde e_i$ and $\bar e$ are the
generators of $H_1(\tilde A\cap\tilde B_i; \Bbb Z)$ and $H_1(A\cap
B; \Bbb Z)$ respectively whose orientations are all induced from $\Bbb S^2$.
It then suffices to show that there is a
neighborhood $V\subset B$ of $p$ such that
$|\rho^{-1}(q)\cap\tilde B_i|\geq k_i$ for any $q\in V\setminus\{p\}$.

For convenience, we let $\tilde A\cap\tilde B_i=\tilde
S^1\times(0,1)$ and $A\cap B=S^1\times(0,1)$. And we can assume that
$$\tilde S^1\times\{t\}\overset{t\to1}\to\longrightarrow\tilde
p_i\text{ and } S^1\times\{t\}\overset{t\to1}\to\longrightarrow
p.\eqno{(10.3)}$$ Let $\mu:\hat S^1\to S^1$ be a $k_i$-covering map
from a circle to a circle, which induces a natural $k_i$-covering
map $\mu:\hat S^1\times(0,1)\to S^1\times (0,1)\ (=A\cap B)$ by
$(u,t)\mapsto (\mu(u),t)$. On the other hand, note that $$H_1(\tilde
A\cap\tilde B_i; \Bbb Z)\cong \pi_1(\tilde A\cap\tilde B_i) \text{
and } H_1(A\cap B; \Bbb Z)\cong \pi_1(A\cap B).$$ Hence, because of
``$\rho_*(\tilde e_i)=k_i\cdot \bar e$'' (see
(10.2)), we know that $\rho|_{\tilde A\cap\tilde B_i}:\tilde
A\cap\tilde B_i\to A\cap B$ has a lifting map $\hat\rho:\tilde
A\cap\tilde B_i\to \hat S^1\times(0,1)$, i.e.
$$\rho|_{\tilde A\cap\tilde B_i}=\mu\circ\hat\rho.$$
Then for any $t$, $\hat\rho(\tilde S^1\times\{t\})$ is homotopy to
$\hat S^1\times(0,1)$. This together with $\rho|_{\tilde A\cap\tilde
B_i}(\tilde S^1\times\{t\})\overset{t\to1}\to\longrightarrow p$ (see
(10.5)) implies that there is a neighborhood $V\subset B$ of $p$
such that $|\rho^{-1}(q)\cap\tilde B_i|\geq k_i$
for any $q\in V\setminus\{p\}$. \hfill$\qed$
\enddemo

\remark{Remark \rm 10.3}  From the proof of Lemma 10.2, one may get a bound on
$\{p\in \Bbb S^2\,|\, |\rho^{-1}(p)|<k\}$, $k\ge 2$. For instance, if there is $p_1\ne p_2\in \Bbb S^2$
such that $\rho^{-1}(p_i)$ is a  single point, then there is no other such kind of points,
i.e.,
$$|\{p\in\Bbb S^2|\ |\rho^{-1}(p)|=1\}|\leq2.$$
It suffices to show that for any $q\in \Bbb
S^2\setminus\{p_1,p_2\}$, $|\rho^{-1}(q)|\geq k$. This can be obtained through the above proof when
one takes $A=\Bbb S^2\setminus\{p_1\}$ (resp. $B=\Bbb S^2\setminus\{p_2\}$)
and $\tilde A=\Bbb S^2\setminus\{\rho^{-1}(p_1)\}$ (resp.
$\tilde B=\Bbb S^2\setminus\{\rho^{-1}(p_2)\}$).
\endremark

\vskip2mm

\demo{Proof of Lemma 10.1}

Claim: {\it For small $t$,
$\pi\circ\Psi_t:\tilde S_c\to S$ is an onto map and there is a full-measure subset $S_c^t$ in
$S$ such that $|(\pi\circ\Psi_t)^{-1}(p)|\geq k$ for any $p\in S_c^t$, where $k$ is the number
in ``$S_k(F_v)$''}. It turns out
that the claim's proof implies that $\tilde Q_c^t\to\tilde Q_c$ as $t\to 0$.

By assuming the claim and that $\tilde Q_c^t\to\tilde Q_c$ as $t\to 0$,
we first prove that $|\tilde Q_c^t|=|\tilde Q_c|$ for small $t$.
Since $\pi\circ\Psi_t$ is distance non-increasing, for $p\in S$
and small $\delta_p>0$, $\bigcup_{\tilde p\in (\pi\circ\Psi_t)^{-1}(p)}B_{\delta_p}(\tilde p)
\subseteq(\pi\circ\Psi_t)^{-1}(B_{\delta_p}(p))$;
and thus the claim implies that
$$\text{Vol}(\tilde S_c) \geq k\cdot\text{Vol} (S),\eqno{(10.4)}$$
and if the equality holds, then
$$\text{Vol}\left((\pi\circ\Psi_t)^{-1}(B_{\delta_p}(p))\right)=
k\cdot\text{Vol} (B_{\delta_p}(p)),\eqno{(10.5)}$$
which implies that $|(\pi\circ\Psi_t)^{-1}(p)|=k$ (resp. $\leq k$) for any $p\in S^t_c$ (resp.
$p\in S\setminus S^t_c$). Meanwhile, we know that the equality in
(10.6) does hold by (8.2.1-2). Then we will use (10.7) to conclude that $|\tilde Q_c^t|=|\tilde Q_c|$ for small $t$.
Since $|\tilde Q_c|<\infty$ (by (8.2.1)) and we have assumed that $\tilde Q_c^t\to\tilde Q_c$ as $t\to 0$,
it suffices to show that each $\tilde q\in \tilde Q_c$ has a neighborhood which
contains only one point of $(\pi\circ\Psi_t)^{-1}(q)\subset \tilde Q_c^t$, where $q=\pi(\tilde q)\in S\setminus S_k(F_v)$.
Since any $\hat q\in (\pi\circ\Psi_t)^{-1}(q)$
converges to some $\tilde q\in \tilde Q_c$ as $t\to0$, from the fact in the proof of (8.1.1) we know that
$\text{diam}(\Sigma_{\hat q}\tilde S_c)\to\pi$ as $t\to 0$ if $\hat q\neq \tilde q$.
This, plus $\dim(S)=\dim(\tilde S_c)=2$ and that $\pi\circ\Psi_t$ is distance non-increasing, implies that if
a small neighborhood of $\tilde q$ contains two (or more) points of
$(\pi\circ\Psi_t)^{-1}(q)$, then for small $\delta>0$
$$\text{Vol}((\pi\circ\Psi_t)^{-1}(B_\delta(q)))>
k\cdot\text{Vol}(B_\delta(q))$$
which contradicts (10.7) (so the proof for $|\tilde Q_c^t|=|\tilde Q_c|$ is finished).

In the rest of the proof, we need only to verify the above claim and $\tilde Q_c^t\to\tilde Q_c$ as $t\to 0$.
Let $S\setminus S_k(F_v)=\{q_1,\cdots,q_m\}$ (see (8.1.1)), and
let $r>0$ be so small that each $B_r(q_j)\subset S$ is a
2-dimensional disc and  $B_r(q_j)\cap B_r(q_{j'})=\emptyset$
for all $j\neq j'$. Then the claim is just the union of the conclusions in the following two cases.

Case 1: $p\in S\setminus\bigcup_{j=1}^mB_\frac r3(q_j)$. By (8.2.2),
we can let $(\pi\circ\Psi_0)^{-1}(p)=\{\tilde
p_1, \cdots, \tilde p_k\}$. And there is
$\hat r$ not depending on $p$ such that $B_{\hat r}(\tilde
p_j)\cap B_{\hat r}(\tilde p_{j'})=\emptyset$ for
all $j\neq j'$, each $B_{\hat r}(\tilde p_j)$ is a 2-dimensional disc, and
$\pi\circ\Psi_0|_{B_{\hat r}(\tilde p_j)}$ is an isometry.
Meanwhile, note that $\pi\circ\Psi_t$ is homotopy to $\pi\circ\Psi_0$.
It then is not hard to conclude that for small $t$ each
$\pi\circ\Psi_t(B_{\hat r}(\tilde p_j))$ contains $p$, and thus
$|(\pi\circ\Psi_t)^{-1}(p)|\geq k$.

Case 2: $p\in \bigcup_{j=1}^mB_\frac r3(q_j)$. Let $q\in\{q_1,\cdots,q_m\}$,
and let $\tilde q\in(\pi\circ\Psi_0)^{-1}(q)$. Note that by (8.2.3)
we have that $\pi\circ\Psi_0(B_r(\tilde q))=B_r(q)$ and
$\pi\circ\Psi_0(\partial B_r(\tilde q))=\partial B_r(q)$
(for all $r$ small). And thus for small $t>0$,
$$\pi\circ\Psi_t\left(B_\frac {3r}4(\tilde q)\right)\subset
B_r(q)\text{ and } \pi\circ\Psi_t\left(B_\frac {3r}4(\tilde q)\setminus
B_\frac r2(\tilde q)\right)\subset
B_r(q)\setminus B_\frac r3(q).$$ Since $\pi\circ\Psi_t$ is homotopy to $\pi\circ\Psi_0$,
the induced homomorphism $(\pi\circ\Psi_t)_*$ from
$$H_i\left(B_\frac {3r}4(\tilde q), B_\frac {3r}4(\tilde q)\setminus
B_\frac r2(\tilde q);\Bbb Z\right)\ (\cong H_i(\Bbb S^2;\Bbb Z))$$ to
$$H_i\left(B_r(q), B_r(q)\setminus B_\frac r3(q);\Bbb
Z\right) \ (\cong H_i(\Bbb S^2;\Bbb Z))$$ satisfies
$$(\pi\circ\Psi_t)_*=(\pi\circ\Psi_0)_*.$$
Moreover, from the end of the proof of (8.1.2), there is an $l\leq k$ such that
$\pi\circ\Psi_0|_{B_r(\tilde q)}:B_r(\tilde q)\to B_r(q)$ is an $l$-branched cover with branched point
$\tilde q$. Hence,
$$(\pi\circ\Psi_t)_*(\tilde e)=(\pi\circ\Psi_0)_*(\tilde e)=l\cdot e,$$
where $\tilde e\in H_2\left(B_\frac {3r}4(\tilde q), B_\frac {3r}4(\tilde q)\setminus
B_\frac r2(\tilde q);\Bbb Z\right)$ and $e\in H_2\left(B_r(q), B_r(q)\setminus B_\frac r3(q);\Bbb
Z\right)$ are generators respectively. Then by Lemma 10.2 and the
arbitrariness of $\tilde q$, it is not hard to conclude that $p\
(\in\bigcup_{j=1}^mB_\frac r3(q_j))$ lies in  $\pi\circ\Psi_t(\tilde
S_c)$ and there is a zero-measure subset $K\subset
\bigcup_{j=1}^mB_\frac r3(q_j)$ such that
$|\pi\circ\Psi_t^{-1}(x)|<k$ for any $x\in K$.

We now show that $\tilde Q_c^t\to\tilde Q_c$ (the union of all $\tilde q$) as $t\to 0$. From the discussions in
Case 1 and 2 above, we can conclude that $\pi\circ\Psi_t$
maps $\tilde S_c\setminus\bigcup_{\text{all }\tilde
q}B_r(\tilde q)$ into $S\setminus\bigcup_{\text{all }
q}B_\frac r2(q)$, and each $B_r(\tilde q)$ contains
at least one point of $(\pi\circ\Psi_t)^{-1}(q)$. It follows that $\Bbb
A_c^t$ (i.e. $\bigcup_{\text{all } q}(\pi\circ\Psi_t)^{-1}(q)$) belongs
to $\bigcup_{\text{all } \tilde q}B_r(\tilde q)$. Note
that we can let $r\to 0$ as $t\to0$, so $\tilde Q_c^t\to\tilde Q_c$ as $t\to 0$.
\hfill$\qed$
\enddemo
We are ready to prove Theorem 0.7.

\demo{Proof of Theorem 0.7}

We first present a proof for $t\in (0,t_0]$ with $t_0$ small.

Let $Q=S\setminus  S_k(F_v)$ and $\tilde Q_c=\varphi_c(Q)$.  By (8.2.2),
$$\pi|_{\tilde
S_c\setminus\tilde Q_c}:\tilde
S_c\setminus\tilde Q_c\to
S\setminus Q\text{ is a covering map}.$$
Observe that
for $\epsilon$ sufficiently small $\pi_1(\tilde S_c\setminus\tilde Q_c)=
\pi_1(\tilde S_c\setminus B_\epsilon(\tilde Q_c))$ because $\tilde Q_c$ consists
of finite points by (8.2.1). And by Lemma 10.1, we may assume $t_0$ so small that for all $t\in [0,t_0]$,
$\tilde Q_c^t\subset B_\epsilon(\tilde Q_c)$ and $|\tilde Q_c^t|=|\tilde Q_c|$.
It follows that
$$\pi_1(\tilde S_c\setminus\tilde Q_c^t)=\pi_1(\tilde S_c\setminus\tilde Q_c)$$
and
$$(\pi\circ\Psi_{t})_*\left(\pi_1(\tilde S_c\setminus\tilde Q_c^t)\right)=(\pi|_{\tilde
S_c})_*\left(\pi_1(\tilde
S_c\setminus\tilde Q_c)\right).$$
Consequently, $\pi\circ\Psi_{t}|_{\tilde
S_c\setminus\tilde Q_c^t}:\tilde S_c\setminus \tilde Q_c^t\to S\setminus Q$
has a lifting map, $$\widetilde{\pi\circ\Psi_{t}}:\tilde
S_c\setminus\tilde Q_c^t\to
\tilde S_c\setminus\tilde Q_c,$$ i.e.,
$\pi\circ\Psi_{t}|_{\tilde S_c\setminus\tilde Q_c^t}=\pi|_{\tilde
S_c}\circ\widetilde{\pi\circ\Psi_{t}}.$ Note that if $c=c_0$, then $\pi|_{\tilde
S_c}=\op{id}$, and thus $\widetilde{\pi\circ\Psi_{t}}=\pi\circ\Psi_{t}.$
Since $\pi\circ\Psi_{t}$ is distance non-increasing and
$\pi|_{\tilde S_c\setminus\tilde Q_c}$ is a local isometry (by (8.2.2)),
$\widetilde{\pi\circ\Psi_{t}}$ is distance non-increasing locally, and thus globally
because $\tilde Q_c^t$ is a finite set (by Lemma
10.1). So, $\widetilde{\pi\circ\Psi_{t}}$ can be
extended uniquely to a distance non-increasing map
$$\widetilde{\pi\circ\Psi_{t}}:\tilde S_c\to\tilde S_c \text{ with }
\pi\circ\Psi_{t}=\pi|_{\tilde
S_c}\circ\widetilde{\pi\circ\Psi_{t}}.$$
Note that $\widetilde{\pi\circ\Psi_0}=\op{id}_{\tilde S_c}$, and that
$\widetilde{\pi\circ\Psi_t}$ is homotopy to $\widetilde{\pi\circ\Psi_0}$
because $\pi\circ\Psi_t$ is homotopy to $\pi\circ\Psi_0$. It then follows that
$\widetilde{\pi\circ\Psi_{t}}$ is an onto map for
$t\in (0,t_0]$. Moreover, we know that $\widetilde{\pi\circ\Psi_{t}}$ is
distance non-increasing, so it is an isometry
for $t\in (0,t_0]$ (ref. [Pet2]). Together with (8.2.3) this implies that,
for all $x\in\Psi_t(\tilde S_c)$ and $r$,
$\pi(B_r(x))=B_r(\pi(x))$.
That is, Theorem 0.7 holds for $t\in (0,t_0]$.
Moreover, as a corollary,
$\pi\circ\Psi_{t_0}:\tilde S_c\to S$ has the same properties as $\pi|_{\tilde S_c}$ in the sense
of Corollary 8.2, where we need to let $\tilde Q_c=\tilde Q_c^{t_0}$.

Next, starting with $\pi\circ\Psi_{t_0}:\tilde S_c\to S$ and
similar to Lemma 10.1, we can obtain that, for $t>t_0$
with $t-t_0$ small, $|\tilde Q_c^t|=|\tilde Q_c^{t_0}|$ and $\tilde Q_c^t\to
\tilde Q_c^{t_0}$ as $t\to t_0$. Similarly, we can verify Theorem 0.7
for $t>t_0$ with $t-t_0$ small. It is clear
that we can give a proof for all $t$ by repeating the process.
\hfill$\qed$
\enddemo

\remark{Remark \rm 10.6}  With some additional work which we intend to carry out elsewhere,
we can  prove (A1) for the case where $\dim(S)=\dim(X)-2$ with
$\dim(X)\geq 5$. From the proof of Theorem A, one sees that the remaining work
is to show that $F_v\setminus S$ divides its a neighborhood into two
components (so we can define $\Psi_{t}$ as (9.1)), and that
$\pi\circ\Psi_{t}$ has a lifting map as in Theorem 0.7.
\endremark

%%%%%%%%%%%%%%%%%%%%%%%%%%%%%%%%%%%%%%%%%%%%%%%%%%%%%%%%%%%%%%%%%%%%%

\vskip4mm

\noindent{\bf Acknowledgement.} The second author would like to thank Professor Hongzhu Gao
for his valuable discussion on the 2-dimensional topology used in the paper.

%%%%%%%%%%%%%%%%%%%%%%%%%%%%%%%%%%%%%%    References    %%%%%%%%%%%%%%%%%%%%%%%%%%

\enddocument